\documentclass[times, review, 10pt]{elsarticle}
\usepackage{amsmath,amssymb,amsfonts}
\usepackage{cases}
\usepackage{color,xcolor}
\usepackage{graphicx}
\usepackage{subfigure}
\usepackage{epstopdf}
\usepackage{multirow}
\usepackage{rotating}
\usepackage{booktabs}
\usepackage{tikz}
\usepackage{color}
\usepackage{makecell}
\usepackage[ruled,linesnumbered]{algorithm2e}
\allowdisplaybreaks[4]

\newtheorem{theorem}{Theorem}[section]
\newtheorem{proposition}{Proposition}[section]
\newtheorem{lemma}{Lemma}[section]

\newtheorem{definition}{Definition}[section]

\graphicspath{{figures/}} 

\newproof{proof}{Proof}

\begin{document}

\begin{frontmatter}

\title{Fast Online $\ell_0$ Elastic Net Subspace Clustering via A Novel Dictionary Update Strategy}

\author[mymainaddress]{Wentao Qu \corref{mycorrespondingauthor}}
\cortext[mycorrespondingauthor]{Corresponding author}
\ead{wtqu96@bjtu.edu.cn}

\author[mymainaddress]{Lingchen Kong}
\ead{lchkong@bjtu.edu.cn}

\author[mysecondaryaddress]{Linglong Kong}
\ead{lkong@ualberta.ca}

\author[mysecondaryaddress]{Bei Jiang}
\ead{bei1@ualberta.ca}

\address[mymainaddress]{School of Mathematics and Statistics, Beijing Jiaotong University, Beijing 100044, P. R. China}
\address[mysecondaryaddress]{Department of Mathematical and Statistical Sciences, University of Alberta, Edmonton T6G 2G1, Canada}

\begin{abstract}
	With the rapid growth of data volume and the increasing demand for real-time analysis, online subspace clustering has emerged as an effective tool for processing dynamic data streams. However, existing online subspace clustering methods often struggle to capture the complex and evolving distribution of such data due to their reliance on rigid dictionary learning mechanisms. In this paper, we propose a novel $\ell_0$ elastic net subspace clustering model by integrating the $\ell_0$ norm and the Frobenius norm, which owns the desirable block diagonal property. To address the challenges posed by the evolving data distributions in online data, we design a fast online alternating direction method of multipliers with an innovative dictionary update strategy based on support points, which are a set of data points to capture the underlying distribution of the data. By selectively updating dictionary atoms according to the support points, the proposed method can dynamically adapt to the evolving data characteristics. Moreover, we rigorously prove the convergence of the algorithm. Finally, extensive numerical experiments demonstrate that the proposed method improves clustering performance and computational efficiency, making it well-suited for real-time and large-scale data processing tasks.
\end{abstract}


\end{frontmatter}


\section{Introduction}\label{1}
The last few decades have witnessed significant advances in data collection technologies, where data typically arrives in real-time and complete datasets are not available at the initial stage. This shift has led to a growing demand for techniques that can efficiently process and adapt to changes in data characteristics \cite{krawczyk2017ensemble, zubarouglu2021data, alam2025online}. 
Online subspace clustering has gained increasing attention for its potential in processing online streaming data \cite{hassani2015subspace,ab2023subspace}. It has been applied in various domains ranging from machine learning to control science, including image processing \cite{rao2012clustering}, video surveillance \cite{hunt2019online}, and real-time anomaly detection \cite{lee2023online}. The core of this technique lies in the ability to capture the subspace structure in the data, which is a key feature of subspace clustering methods.

Subspace clustering is a classic clustering method to deal with high-dimensional data. It aims to identify the underlying subspace in the data space and use this information to derive the clustering results \cite{abdolali2021beyond,qu2023survey}. 
Due to the ability to handle arbitrary shape data and clear mathematical principles, spectral-based methods have become predominant in subspace clustering \cite{chen2020stochastic}. It is crucial to construct a high-quality representation matrix to capture the similarity between samples in spectral-based subspace clustering methods. Existing studies have focused on leveraging the self-expression property of data to get the representation matrix, achieving promising results. So far, various subspace clustering methods have been developed by utilizing the prior information of different data structures. 
Elhamifar and Vidal \cite{Elhamifar2013} introduced the sparse compressed sensing techniques into subspace clustering, which utilizes the sparsity of the representation matrix to reveal the local information of data. 
In contrast, to capture the global structure of data, Liu \textit{et al.} \cite{Liu2013lrr} presented a low-rank representation subspace clustering (LRSC) method based on the low-rank regularization. 
Lu \textit{et al.} \cite{Lu2012} developed a subspace clustering method via least squares regression, and proved that their method can achieve the grouping effect by leveraging the correlation of data. 
You \textit{et al.} \cite{you2016oracle} proposed a subspace clustering via elastic net regularization, which can balance the subspace preservation and connectivity. Feng \textit{et al.} \cite{feng2024adaptive} presented an adaptive weighted dictionary representation using anchor graph for subspace clustering. Most recently, Lu \textit{et al.} \cite{Lu2019bdr} pointed out that a good subspace representation matrix should have a block diagonal structure. Base on this, Lin and Song \cite{lin2023convex} proposed a convex subspace clustering by adaptive block diagonal representation. 

Although the aforementioned methods have achieved promising performance, they all use the data itself as the dictionary, which significantly limits their application to high-dimensional data. This limitation has become increasingly evident as the magnitude of data continues to grow \cite{paul2020online}. 
In particular, many modern applications involve online dynamic data, where the entire dataset is not available upfront but arrives incrementally.  Consequently, batch subspace clustering methods are not suitable for these scenarios. This has driven the development of online subspace clustering \cite{shen2016online}, a variant of subspace clustering specifically designed to efficiently handle online dynamic data. 
The main challenges faced by online subspace clustering are as follows \cite{bifet2023machine}: The memory is limited as it cannot store the entire dataset at once; They need to react quickly to incoming data; They should be able to detect and handle data corruption to ensure robustness. To address these challenges, existing online subspace clustering methods have primarily focused on developing effective dictionary update rules. 
One approach is to combine subspace clustering with dictionary learning. 
For instance, Shen \textit{et al.} \cite{shen2016online} extended the LRSC by introducing a nonconvex reconstruction of the nuclear norm and combining it with dictionary learning. 
Madden \textit{et al.} \cite{madden2019online} adapted the sparse subspace clustering into an online framework without revisiting the whole dataset. 
Chen \textit{et al.} \cite{chen2023two} integrated the $\ell_0$ sparse subspace clustering and dictionary learning into a unified framework, and provided a dictionary initialization strategy for improved performance. 
For the tensor online streaming data, Li \textit{et al.} \cite{li2019online} proposed a low-rank tensor subspace clustering based on dictionary learning, and designed a stochastic optimization algorithm to handle dynamic data. Based on this, Wu \cite{wu2023online} incorporated tensor nuclear norm decomposition into the above framework for enhanced robustness. 
Beyond dictionary learning, the challenges posed by online data can also be addressed through a two-stage process. For example, Lee \textit{et al.} \cite{lee2016online} reconstructed the LRSC by the nonconvex reconstruction of the nuclear norm and developed a dictionary update rule based on sparsification and pruning procedures. 
Li \textit{et al.} \cite{li2018online} proposed a two-stage algorithm for LRSC, which first extracted the subspace structure from a small subset of data and then utilizes this structure to compute the intrinsic principal components of the entire dataset. 
Sui \textit{et al.} \cite{sui2022dynamic} proposed a dynamic low-rank subspace clustering method, which contains static learning to select the dictionary to adapt to the data distribution. 
Chen \textit{et al.} \cite{chen2024online} further advanced this line of research by incorporating a low-dimensional projection into the $\ell_0$ sparse representation clustering and designed a dictionary sample selection scheme to efficiently update the dictionary.  
Although the existing online subspace clustering methods have shown some success in practice, they still struggle to effectively capture the complex and evolving structures present in dynamically changing online data streams. 
More importantly, the current dictionary update rules for online subspace clustering fail to adequately reflect the true distribution of the data, limiting their effectiveness in representing the underlying structure of incoming data. To address this limitation, the support points, an important research topic in statistics, can be employed. 
Mak and Joseph \cite{mak2018support} proposed a new method for data compressing using support points, obtained by minimizing the energy distance. Recently, the concept of support points has gained popularity in data sampling, see \cite{joseph2022split,song2024large}. This offers us a new strategy for building the dictionary in subspace clustering, enabling a better representation of the underlying data structure.

Motivated by the above observation, in this paper, we propose a novel subspace clustering model by integrating the $\ell_{0}$ norm and Frobenius norm into a unified framework, called $\ell_0$-ENSC. The $\ell_0$-ENSC model not only performs variable selection, but also exhibits the block diagonal property, meaning it can effectively group highly correlated samples. Our model is an NP-hard nonconvex discontinuous problem due to the discrete nature of the $\ell_0$ norm penalty term. 
We rigorously establish the existence of optimal solutions and derive optimality conditions for the $\ell_0$-ENSC model. 
Based on these theoretical insights, we design an efficient online alternating direction method of multipliers (ADMM). Furthermore, to better handle online data, we incorporate a dictionary updating strategy based on support points into the ADMM framework. 
Finally, extensive numerical experiments conducted on six public datasets demonstrate the superior performance of our method in comparison with several state-of-the-art online subspace clustering methods. 

The contributions of this paper are summarized as follows.
\begin{enumerate}[(1)]
	\item  It proposes an $\ell_0$ elastic net subspace clustering model ($\ell_0$-ENSC), and proves that the proposed model has the block diagonal property. It also establishes the relationship between the Karush-Kuhn-Tucker (KKT) point, P-stationary point, and local optimal solutions. This provides theoretical support for algorithm design.
	
	\item It develops a fast online algorithm that integrates a dictionary update strategy based on support points into the ADMM framework to dynamically capture online data distributions. Theoretically, it is rigorously proved that the proposed algorithm converges to the P-stationary point, which is also a uniquely local minimizer. 
	 
	\item It demonstrates the superior performance of the $\ell_0$-ENSC method compared to the state-of-the-art online subspace clustering methods on six real-world datasets. Our method improves both clustering performance and computational efficiency, highlighting its advantages in practical applications.
\end{enumerate}

The remainder of this paper is organized as follows. 
The proposed $\ell_0$-ENSC model is presented in Section \ref{sec_mod}. Section \ref{sec_alg} designs an online subspace clustering algorithm with a novel dictionary update strategy. The extensive experimental results are discussed in Section \ref{sec_exp}. Finally, Section \ref{sec_con} concludes this paper and proposes some potential future research directions.

We end this section with the notations used throughout this paper. All matrices are represented by uppercase letters (e.g., $X$), and all vectors are represented by lowercase bold letters (e.g., $\mathbf{x}$). For a vector $\mathbf{x}\in\mathbb{R}^m$, ${x}_i$ denotes its $i$-th element, and $\|\mathbf{x}\|_2$ represents the Euclidean norm of $\mathbf{x}$. 
For a matrix $X\in\mathbb{R}^{m\times n}$, we denote the $i$-th column by $\mathbf{x}_i \in \mathbb{R}^{m}$ and the $(i, j)$-th element by $X_{ij} \in \mathbb{R}$. $\|X\|_{F}$ denotes the Frobenius norm of $X$, while  $\|X\|_{0}$ represents the $\ell_{0}$ (pseudo) norm which counts the number of nonzero elements of $X$. 
Additional notation will be introduced wherever it appears.

\section{$\ell_0$ Elastic Net Subspace Clustering}\label{sec_mod}
In this section, we begin by introducing the proposed $\ell_0$-ENSC model. Then we demonstrate the block diagonal property of the $\ell_0$-ENSC model. Finally, we establish the optimality conditions of the proposed model.

\subsection{Model Formulation}
To effectively capture both local and global structure in high-dimensional data, we establish the $\ell_0$ elastic net subspace clustering ($\ell_0$-ENSC) model. The $\ell_0$-ENSC model aims to leverage the strengths of both $\ell_0$ sparse representation and group effect regularization to provide a more robust solution for subspace clustering, especially in scenarios involving high-dimensional datasets with complex structures. Specifically, we consider solving the following problem
\begin{eqnarray}\label{OENSC}
	\min_{X\in \mathbb{R}^{m\times n}}~~\frac{1}{2}\lVert Z-DX \rVert _{F}^{2}+\lambda _1\lVert X \rVert _0+\frac{\lambda _2}{2}\lVert X \rVert _{F}^{2},
\end{eqnarray}
where $Z\in \mathbb{R}^{p\times n}$ represents the dataset with $n$ samples, $D\in \mathbb{R}^{p\times m}$ is the dictionary, and $X\in \mathbb{R}^{m\times n}$ denotes the corresponding coefficient matrix. $\lambda _1,\lambda _2\geq 0$ are the tuning parameters that control the trade-off between the regularization terms.  
In the objective function of problem \eqref{OENSC}, the first term represents the subspace self-representation property, implying that each point can be represented as a linear combination of others through the dictionary $D$. The $\ell_0$ term is used to promote the local sparsity of $X$, thus enhancing model interpretability and reducing model complexity. This sparse representation allows the model to effectively capture the local structure within the subspaces. The third term encourages the global group effect, which tends to shrink coefficients of correlated data and group them together. Thus, our model can capture the underlying relationships between different subspaces and improve clustering quality. This combined regularization makes the proposed model well-suited for efficiently clustering data with complex relationships while maintaining both sparsity and robustness. 

Notably, several existing models can be viewed as a special case of our model. For example, EnSC \cite{you2016oracle} can be interpreted as a convex relaxation of our model by replacing the $\ell_0$ regularization with the $\ell_1$ regularization. The proposed $\ell_0$-ENSC \eqref{OENSC} reduces to the LSR model \cite{Lu2012} when $\lambda_1 = 0$, and degenerates into the $\ell_0$-SSC model \cite{Elhamifar2013} when $\lambda_2 = 0$.

Theorem \ref{th_bdq} further elaborates on the property of our $\ell_0$-ENSC that its solution exhibits the block diagonal property. The detailed proof is provided in the Appendix.

\begin{theorem}\label{th_bdq}
	Assume the dataset $Z =[Z_1, \ldots, Z_c]\in \mathbb{R}^{p\times n}$ drawn from $c$ independent subspaces $\{\mathcal{S}_i\}_{i = 1}^c$, where $Z_i\in\mathbb{R}^{p\times n_i}$ represents data from $\mathcal{S}_i$ with $\sum_{i=1}^c n_i = n$. 
	Let $\Omega$ be a set consisting of square matrices.
	If $X^*\in \Omega$ is an optimal solution of \eqref{OENSC}, then $X^*$ satisfies the block diagonal property, that is
	\begin{eqnarray}\label{BQ}
		X^* = \left( \begin{matrix}
			X_1^* & 0 & \cdots & 0 \\ 0 & X_2^* & \cdots & 0\\ \vdots & \vdots & \ddots & \vdots \\ 0 & 0 & \cdots & X_c^*
		\end{matrix}\right),
	\end{eqnarray}
	with $X_i\in\mathbb{R}^{m_i\times n_i}$ corresponding to $Z_i$.
\end{theorem}

\subsection{Optimality Conditions}
In this section, we first define the Karush-Kuhn-Tucker (KKT) point and P-stationary point of the $\ell_0$-ENSC model. Then the relationship between the KKT points, P-stationary point, and the local optimal solution is established.

To better handle the nonconvex term, we introduce an auxiliary variable $Y$ and then equivalently reformulate the $\ell_0$-ENSC model \eqref{OENSC} into the following constrained optimization problem
\begin{eqnarray}\label{OENSC1}
	\begin{aligned}
		\min_{Y, X }&~~ f(Y,X) := \frac{1}{2}\lVert Z-DX \rVert _{F}^{2}+\lambda _1\lVert Y \rVert _0+\frac{\lambda _2}{2}\lVert X \rVert _{F}^{2},\\
		\text{s.t.} & ~~Y = X.
	\end{aligned}
\end{eqnarray}
For convenience, denote $h(X):= \frac{1}{2}\lVert Z-DX \rVert _{F}^{2}+\frac{\lambda_2}{2}\lVert X \rVert _{F}^{2}$. Then, $\nabla h(X) = (D^{\top}D+\lambda_2 I) X-D^{\top}Z$. It is easy to verify that $h$ is strongly convex and $r$-smooth with $r = \lVert D^{\top}D\rVert_2 + \lambda_2$. For a $r$-smooth function $h$, the following inequality holds
\begin{eqnarray}\label{eq_smooth}
	h(X)\leq h(Y) + \left\langle \nabla h(Y), X -Y\right\rangle +\frac{1}{r}\|X -Y\|_F^2.
\end{eqnarray}
Then the Lagrangian function of \eqref{OENSC1} is defined as
\begin{eqnarray}
	\mathcal{L}(Y,X,U) = h(X) +\lambda _1\lVert Y \rVert _0 + \langle U, X-Y\rangle ,
\end{eqnarray}
where $U=[\mathbf{u}_1, \ldots, \mathbf{u}_n]\in\mathbb{R}^{m\times n}$ is the Lagrangian multiplier. We say that $(X^*,Y^*, U^*)$ is a KKT point of model \eqref{OENSC1} if it satisfies 
\begin{eqnarray}\label{KKT}
	\begin{aligned}
		\left\{ \begin{aligned}			
			& 0 \in \lambda_1 \partial\lVert Y^* \rVert _0-U^*, \\ 
			& 0 = \nabla h(X^*) + U^*,\\
			& Y^* = X^*,\\ 
		\end{aligned} \right.
	\end{aligned}
\end{eqnarray}
where $\partial\lVert Y^* \rVert _0 = \{V : V_{ij}=0, \forall (i,j)\in \mathcal{T}(Y^*)\}$ with $\mathcal{T}(Y)$ being the support set of $Y$. Note that it is difficult to find a KKT point of \eqref{OENSC1}. Thus, we give the definition of the P-stationary point, which can be easily calculated by the proximal operator.
\begin{definition}
	For a given $\alpha>0$, a point $(Y^*, X^*, U^*)$ is called a P-stationary point of problem \eqref{OENSC1} if it satisfies
	\begin{eqnarray}\label{stationary}
		\begin{aligned}
			\left\{ \begin{aligned}
				& Y^* = \operatorname{Prox}_{\alpha\lambda_1\|\cdot\|_0}\left(Y^*+\alpha U^*\right)\\
				& 0 =  \nabla h(X^*) + U^*, \\
				& Y^* = X^*,\\  
			\end{aligned} \right.
		\end{aligned}
	\end{eqnarray}
	where $ \operatorname{Prox}_{\alpha\lambda_1\|\cdot\|_0}\left(\cdot\right)$ is the proximal operator of $\alpha\lambda_1\|\cdot\|_0$.
\end{definition}

It is worth mentioning that  $\operatorname{Prox}_{\alpha\lambda_1\| \cdot \|_0}\left( \cdot \right)$ is a multi-valued operator due to the nonconvexity of $\alpha\lambda_1\|\cdot\|_0$. According to the definition of the proximal operator, for a given matrix $A$, we know that, for any $i\in\{1, 2, \ldots, n\}$, 
\begin{eqnarray}
	\left[\operatorname{Prox}_{\alpha\lambda_1\|\cdot\|_0}\left(A\right)\right]_{ij} = \underset{X_{ij}}{\arg\min}  \frac{1}{2}\|X_{ij}-A_{ij} \|_2^2 + \alpha\lambda_1 \mathbb{I} (X_{ij}\neq 0) . 
\end{eqnarray}
Without loss of generality, we can calculate the proximal operator of $\|\cdot\|_0$ as
\begin{equation}\label{eq_prox0}
	\left[\operatorname{Prox}_{\alpha\lambda_1\|\cdot\|_0}\left(A\right)\right]_{ij} = \begin{cases}
		0, & \frac{1}{2}A_{ij}^2\leq \alpha\lambda_1,\\
		A_{ij}, & \frac{1}{2}A_{ij}^2> \alpha\lambda_1.\\
	\end{cases}
\end{equation}

The following proposition provides a more explicit formulation of the stationary point defined by \eqref{stationary}. Detailed proofs of all theoretical results are given in the Appendix.

\begin{proposition}\label{pop_1}
	For a given $\alpha>0$, if $(Y^*,X^*, U^*)$ is a P-stationary point of problem \eqref{OENSC1}. Then the following two assertions hold
	\begin{enumerate}[(a)]
		\item $\left| Y^*_{ij} \right| > \sqrt{2\alpha\lambda_1}$ or $Y^*_{ij} = 0, \forall i, j$;
		\item $\left| X^*_{ij} \right| > \sqrt{2\alpha\lambda_1}$ or $X^*_{ij} = 0, \forall i, j$.
	\end{enumerate}
\end{proposition}

To conclude this section, we establish the relationships between the P-stationary point, KKT point, and local optimal solutions of problem \eqref{OENSC1}. The subsequent theorem provides the first-order necessary condition of \eqref{OENSC1}.  
\begin{theorem}\label{th_nc}
	Suppose that $(Y^*,X^*)$ is a local minimizer of \eqref{OENSC1}. Then, there exists $U^*$ such that $(Y^*,X^*, U^*)$ is a KKT point. 
	Furthermore, there exists an $\epsilon_* > 0$ such that $(Y^*,X^*, U^*)$ is a P-stationary point with $0<\alpha< \alpha_*:=\min \{\alpha_1^*,\alpha_2^*\}$, where 
	\begin{eqnarray*}
		\alpha_1^* :=\frac{\epsilon_* }{2\lambda _1\lVert Y^* \rVert _0 + 2\sqrt{\lambda _1^2\lVert Y^* \rVert _0^2 +  \epsilon_* \lVert\nabla h(X^*)\rVert _F^2}}, \alpha_2^* :=\frac{1}{\sqrt{n}\|D^{\top}D  + \lambda_2 I\|_F}.
	\end{eqnarray*}
\end{theorem}

Building upon the result above, the next theorem shows that a P-stationary point of \eqref{OENSC1} must be a KKT point and also a uniquely local minimizer.

\begin{theorem}\label{th_skl}
	For a given $\alpha>0$, if $(Y^*,X^*, U^*)$ is a P-stationary point of problem \eqref{OENSC1}. Then it is a KKT point and also a uniquely local minimizer of \eqref{OENSC1}.
\end{theorem}

Theorem \ref{th_skl} asserts that the P-stationary point is a sharper optimal condition than the KKT point. Most importantly, the optimality of the P-stationary point can be guaranteed and easy to verify.

\section{Online Optimization Algorithm}\label{sec_alg}
This section first presents an online $\ell_0$ elastic net subspace clustering algorithm based on ADMM. Then a strategy for updating the dictionary based on support points is introduced. Finally, the dictionary update strategy is incorporated into the online algorithm to improve overall efficiency and performance.

\subsection{Online Optimization Algorithm}
In this part, we employ the famous alternating direction method of multipliers (ADMM) to solve model \eqref{OENSC1}, which contains three steps in each iteration. The first step is to update $X$ by solving a quadratic programming. The second step is to update $Y$ by computing a hard threshold operator, which induces the sparsity of $Y$. The correction of the Lagrange multiplier is given in the last step.

Denote $Z = \{\mathbf{z}_1, \ldots, \mathbf{z}_n\}$, $X = \{\mathbf{x}_1, \ldots, \mathbf{x}_n\}$, and $Y = \{\mathbf{y}_1, \ldots, \mathbf{y}_n\}$. Then, the problem \eqref{OENSC1} is now separable into different time slots and can be rewritten as
\begin{eqnarray}\label{OENSC_ve}
	\begin{aligned}
		\min_{\{\mathbf{y}_i\}, \{\mathbf{x}_i\}}&~~ \frac{1}{n}\sum_{i=1}^{n} \left(\frac{1}{2}\lVert \mathbf{z}_i-D\mathbf{x}_i \rVert _{2}^{2}+\lambda _1\lVert \mathbf{y}_i \rVert _0 +\frac{\lambda _2}{2}\lVert \mathbf{x}_i \rVert _{2}^{2}\right),\\
		\text{s.t.} & ~~\mathbf{y}_i=\mathbf{x}_i, \forall i.
	\end{aligned}
\end{eqnarray}
The augmented Lagrangian function is defined as

\begin{eqnarray}\label{lag}
	\mathcal{L}_{\sigma}(Y, X, U) = \frac{1}{n}\sum_{i=1}^{n} \mathcal{L}_{\sigma}(\mathbf{y}_i,\mathbf{x}_i,\mathbf{u}_i),
\end{eqnarray}
where $\mathcal{L}_{\sigma}(\mathbf{y}_i,\mathbf{x}_i, \mathbf{u}_i)=\frac{1}{2}\lVert \mathbf{z}_i-D\mathbf{x}_i \rVert _{2}^{2}+\lambda _1\lVert \mathbf{y}_i \rVert _0 +\frac{\lambda _2}{2}\lVert \mathbf{x}_i \rVert _{2}^{2}+ \langle \mathbf{u}_i, \mathbf{x}_i-\mathbf{y}_i\rangle + \frac{\sigma}{2}\lVert \mathbf{x}_i-\mathbf{y}_i \rVert _{2}^{2}$ and $\sigma> 0$ is a given parameter. Given the $k$-th iteration $(\mathbf{x}_i^k,\mathbf{y}_i^k, \mathbf{z}_i^k )$, the iteration schemes of ADMM for solving the optimization problem \eqref{OENSC1} are given by
\begin{numcases}{}
	\mathbf{y}_i^{k+1} = \underset{\mathbf{y}_i}{\arg\min} \mathcal{L}_{\sigma}(\mathbf{y}_i,\mathbf{x}_i,\mathbf{u}_i),\label{eq:y1}\\
	\mathbf{x}_i^{k+1} = \underset{\mathbf{x}_i}{\arg\min} \mathcal{L}_{\sigma}(\mathbf{y}_i,\mathbf{x}_i,\mathbf{u}_i), \label{eq:x1}\\ 
	\mathbf{u}_i^{k+1} =\mathbf{u}_i^k + \sigma (\mathbf{x}_i^{k+1} - \mathbf{y}_i^{k+1}).\label{eq:u} 
\end{numcases}
Below we show how to compute $\mathbf{x}_i^{k+1}$ and $\mathbf{y}_i^{k+1}$ in closed-form solutions.

\textbf{(i) Update $\mathbf{y}_i^{k+1}$.} When $\{\mathbf{x}_i\}$ and $\mathbf{u}_i$ are fixed, the solution of \eqref{eq:y1} can be obtained  as follows
\begin{eqnarray}
	\begin{aligned}
		\mathbf{y}_i^{k+1} &= \underset{\mathbf{y}_i}{\arg\min} \left\{ \lambda _1\lVert \mathbf{y}_i \rVert _0 + \frac{\sigma}{2}\lVert \mathbf{x}_i^k + \mathbf{u}_i^k/\sigma  -\mathbf{y}_i \rVert _{2}^{2}\right\}\\
		& = \operatorname{Prox}_{\frac{\lambda_1}{\sigma}\|\cdot\|_0}(\mathbf{x}_i^k + \mathbf{u}_i^k/\sigma)
	\end{aligned}
\end{eqnarray}
By using the hard threshold operator, we can derive the closed-form solution of $Y$ as
\begin{equation}\label{eq:y}
	\mathbf{y}_i^{k+1} = \mathcal{H}_{\sqrt{\frac{2\lambda_1}{\sigma}}}\left(\mathbf{x}_i^k+\frac{\mathbf{u}_i^k}{\sigma}\right),
\end{equation}
where $\mathcal{H}_{t}\left(x\right) = \begin{cases}
	0, & |x|\leq t,\\
	x, & |x|> t.
\end{cases}$.

\textbf{(ii) Update $\mathbf{x}_i^{k+1}$.} When $\mathbf{y}_i$ and $\mathbf{u}_i$ are fixed, the $\mathbf{x}_i$-subproblem \eqref{eq:x1} is rewritten as
\begin{eqnarray}\label{eq:x0}
	\min_{\mathbf{x}_i} \left\{\frac{1}{2}\lVert \mathbf{z}_i-D\mathbf{x}_i \rVert _{2}^{2} +\frac{\lambda _2}{2}\lVert \mathbf{x}_i \rVert _{2}^{2}+ \langle \mathbf{u}_i^{k}, \mathbf{x}_i\rangle + \frac{\sigma}{2}\lVert \mathbf{x}_i-\mathbf{y}_i^{k+1} \rVert _{2}^{2}\right\}.
\end{eqnarray}
Since this is a strongly convex function, we can easy get its closed-form solution as
\begin{eqnarray}\label{eq:x}
	\mathbf{x}_i^{k+1} = \left[ D^{\top}D+(\lambda_2+\sigma)I\right]^{-1}\left( D^{\top}\mathbf{z}_i-\mathbf{u}_i^{k}+\sigma \mathbf{y}_i^{k+1} \right).
\end{eqnarray}

Based on the above analysis, the iteration framework of ADMM for solving \eqref{OENSC1} is summarized in Algorithm \ref{alg1}.

\begin{algorithm}[t]
	\caption{Online ADMM for Solving \eqref{OENSC1}} 
	\label{alg1}
	\KwIn{ $Z\in\mathbb{R}^{p\times n}, D\in\mathbb{R}^{p\times m}$.}
	\textbf{Initialize:} $X^0, Y^0, U^0\in \mathbb{R}^{m\times n}, \lambda_1, \lambda_2$, and $\sigma$.\\
	\For{$i=1,2,\ldots $}{
		\While{no converged}{ 
			Update $\mathbf{y}_i^{k+1}$ by \eqref{eq:y};\\
			Update $\mathbf{x}_i^{k+1}$ by \eqref{eq:x};\\
			Update $\mathbf{u}_i^{k+1}$ by \eqref{eq:u}.
		}
	}
\end{algorithm}

Next, we shall establish the convergence of Algorithm \ref{alg1}. Before doing this, we present two useful lemmas. For convenience, denote $\mathbf{w}_i = \{\mathbf{y}_i, \mathbf{x}_i, \mathbf{u}_i\}$. 

\begin{lemma}[Sufficient Decrease Lemma]\label{nonincreasing}
	Let $\{\mathbf{w}_i^k\}$ be the sequence  generated by Algorithm \ref{alg1}. Then the generated augmented Lagrangian sequence is nonincreasing, i.e.,
	\begin{eqnarray}
		\mathcal{L}_{\sigma}(\mathbf{w}_i^{k+1})- \mathcal{L}_{\sigma}(\mathbf{w}_i^{k}) \leq -\kappa\|\mathbf{x}_i^{k+1}-\mathbf{x}_i^k\|_2^2, 		
	\end{eqnarray}
	where $\kappa = \frac{(2+\sigma)(\gamma+\lambda_2)+\sigma^2}{2\sigma}$ and $\gamma = \lambda_{min}(D^{\top}D)$.
\end{lemma}
\begin{lemma}\label{bound}
	Let $\{\mathbf{w}_i^k\}$ be the sequence generated by Algorithm \ref{alg1} and $\sigma \geq  r$. Then,  $\{\mathbf{w}_i^k\}$ is bounded. Moreover,  the following statement holds
	\begin{eqnarray}
		\lim\limits_{k\rightarrow \infty} \|\mathbf{w}_i^{k+1}-\mathbf{w}_i^k\|_2=0.
	\end{eqnarray}
\end{lemma}

Lemma \ref{bound} shows that both the generated sequence $\{\mathbf{w}_i^k\}$ and the augmented Lagrangian sequence $\{\mathcal{L}_{ \sigma}(\mathbf{w}_i^k)\}$ are bounded. Combining with the nonincreasing of $\{\mathcal{L}_{ \sigma}(\mathbf{w}_i^k)\}$ in Lemma \ref{nonincreasing}, it follows that $\{\mathcal{L}_{ \sigma}(\mathbf{w}_i^k)\}$ is convergent. With the above two lemmas, we are ready to establish the main convergence result in the following theorem.

\begin{theorem}\label{th_con}
	Let $\{\mathbf{w}_i^k\}$ be the sequence generated by Algorithm \ref{alg1}  and $\sigma \geq  r$. Then any accumulation point of $\{\mathbf{w}_i^k\}$ is a P-stationary point and also a uniquely local minimizer of \eqref{OENSC1}.
\end{theorem}

Next, we analyze the complexity of Algorithm \ref{alg1}. Since Algorithm \ref{alg1} is an online method that processes samples in a streaming manner, its computational complexity per sample is given by $\mathcal{O}\left( m(p^2 + m^2) \right)$, where $m$ and $p$ are the numbers of dictionary samples and features respectively. This complexity arises from the following three main components. 
Firstly, updating $\mathbf{y}_i$ primarily involves element-wise addition and a hard threshold operation, resulting in a computational complexity of $\mathcal{O}\left( m \right)$. 
Secondly, updating $\mathbf{x}_i$ is the most computationally intensive step. This step involves solving an equation that requires matrix operations such as calculating $D^{\top}D$, adding a scaled identity matrix, and computing the inverse of the resulting matrix. Hence, the overall complexity for updating $\mathbf{x}_i$ is $\mathcal{O}\left( m(p^2 + m^2) \right)$. 
Finally, updating $\mathbf{u}_i$ only involves simple vector addition, which has a complexity of $\mathcal{O}\left(m\right)$. 

\subsection{Dictionary Update Based on Support Points}
The dynamic nature of online data requires that the model update its parameters efficiently without reprocessing the entire dataset. To meet this need, a novel dictionary update strategy based on support points is designed. This dictionary update strategy not only enhances the adaptability of the model by selectively updating the dictionary atoms that best represent the current data characteristics but also significantly enhances the computational efficiency.

We begin by recalling the definition of the support points proposed in \cite{mak2018support}.

\begin{definition}[Support Points]
	For a dataset $Z\in \mathbb{R}^{p\times n}$, the support points of $Z$ are defined as
	\begin{eqnarray}\label{sp}
		\{\boldsymbol{\xi}_i\}_{i=1}^{m}  = \underset{\mathbf{d}_1, \ldots, \mathbf{d}_m}{\arg\min} \frac{2}{mn}\sum_{i=1}^{m}\sum_{l=1}^{n} \|\mathbf{z}_l-\mathbf{d}_i\|_2- \frac{1}{m^2}\sum_{i=1}^{m}\sum_{j=1}^{m} \|\mathbf{d}_i-\mathbf{d}_j\|_2.
	\end{eqnarray}
\end{definition}
It follows from the above definition that the support points are intended to reconstruct the dataset $Z$ which best represents $Z$ with respect to goodness-of-fit. Moreover, the quality of reconstruction depends on the number of support points ($m$). 

Based on the geometric fact that a $d$-dimensional space requires at least $d + 1$ points to span it, we derive the following proposition.
\begin{proposition}\label{th_sp_m}
	Assume the dataset $Z\in \mathbb{R}^{p\times n}$ is drawn from $c$ independent subspaces $\{\mathcal{S}_i\}_{i = 1}^c$ with the corresponding dimension $\{ d_i = \operatorname{dim}(\mathcal{S}_i) \}_{i = 1}^c$. Then, it needs at least $\sum_{i=1}^c \left( d_i + 1 \right)$ points to adequately capture the distribution of the entire data.
\end{proposition}

In practice, the potential subspaces are usually disjoint, that is the independence condition in Proposition \ref{th_sp_m} may not hold. Consequently, for a specific application, the required number of support points is typically much less than $\sum_{i=1}^c \left( d_i + 1 \right)$.

It is worth pointing out that the objective function of \eqref{sp} is a difference of convex functions. Hence, we use the convex-concave procedure (CCP) to solve \eqref{sp}. The main idea of CCP is to first substitute  the concave component in the difference of convex objective with a convex upper bound, then solve the resulting ``surrogate" formulation (which is convex) by convex programming techniques.

Now we turn to derive the closed-form of \eqref{sp}. Firstly, we majorize the concave component $-{m^{-2}}\sum_{i=1}^{m}\sum_{j=1}^{m} \|\mathbf{d}_i-\mathbf{d}_j\|_2$ by first-order Taylor expansion at the current iterate $\{\mathbf{d}^{'}_i\}_{i=1}^l$, yielding the following surrogate convex program
\begin{eqnarray}\label{scp}
		\min_{\mathbf{d}_1, \ldots, \mathbf{d}_m}\frac{2}{mn}\sum_{i=1}^{m}\sum_{l=1}^{n} \|\mathbf{z}_l-\mathbf{d}_i\|_2 - \frac{1}{m^2}\left[ \sum_{i=1}^{m}\sum_{j=1}^{m} \left( \|\mathbf{d}^{'}_i-\mathbf{d}^{'}_j\|_2 + \frac{2(\mathbf{d}_i-\mathbf{d}^{'}_j)^{\top}(\mathbf{d}^{'}_i-\mathbf{d}^{'}_j)}{\|\mathbf{d}^{'}_i-\mathbf{d}^{'}_j\|_2}  \right)\right].
\end{eqnarray}

Note that, the first term in \eqref{scp} is a nonsmooth term, which greatly limits the ability to solve the problem efficiently. To deal with this issue, the following lemma gives us an idea to further convexify the first term in \eqref{scp}.
\begin{lemma}[Convexification]\label{le_convex}
	$Q(\mathbf{d}|\mathbf{d}^{'}) = \frac{\|\mathbf{d}\|_2^2}{2\|\mathbf{d}^{'}\|_2}+\frac{\|\mathbf{d}^{'}\|_2}{2}$ majorizes $\|\mathbf{d}\|_2$ at $\mathbf{d}^{'}$ for any $\mathbf{d}^{'}\in\mathbb{R}^{M}$.
\end{lemma}

By the above lemma, problem \eqref{scp} can be transformed into
\begin{eqnarray}
	\begin{split}
		\min_{\mathbf{d}_1, \ldots, \mathbf{d}_m}\frac{2}{mn}\sum_{i=1}^{m}&\sum_{l=1}^{n} \left[ \frac{\|\mathbf{z}_l-\mathbf{d}_i\|_2^2}{2\|\mathbf{z}_l-\mathbf{d}^{'}_i\|_2}+\frac{\|\mathbf{z}_l-\mathbf{d}^{'}_i\|_2}{2} \right] \\
		& - \frac{1}{m^2}\left[ \sum_{i=1}^{m}\sum_{j=1}^{m} \left( \|\mathbf{d}^{'}_i-\mathbf{d}^{'}_j\|_2 + \frac{2(\mathbf{d}_i-\mathbf{d}^{'}_j)^{\top}(\mathbf{d}^{'}_i-\mathbf{d}^{'}_j)}{\|\mathbf{d}^{'}_i-\mathbf{d}^{'}_j\|_2}  \right)\right].
	\end{split}
\end{eqnarray}
By taking the derivative of the above equation with respect to $\mathbf{d}_i$ and setting it equal to zero, the global minimizer of $D$ can be obtained directly as
\begin{eqnarray}\label{eq:d}
	\mathbf{d}_i = q^{-1}(\mathbf{d}_i^{'};\{\mathbf{z}_l\}_{l=1}^n) \times \left(\frac{n}{m}\sum_{j=1,j\neq i}^{m}\frac{\mathbf{d}_i^{'}-\mathbf{d}_j^{'}}{\|\mathbf{d}_i^{'}-\mathbf{d}_j^{'}\|_2} + \sum_{l=1}^{n}\frac{\mathbf{z}_l}{\|\mathbf{d}_i^{'}-\mathbf{z}_l\|_2}\right),
\end{eqnarray}
where $q(\mathbf{d}_i^{'};\{\mathbf{z}_l\}_{l=1}^n)=\sum_{l=1}^{n}\|\mathbf{d}_i^{'}-\mathbf{z}_l\|_2^{-1}$ with $i=1,\ldots, m$.

In summary, the whole framework of solving \eqref{sp} is outlined in Algorithm \ref{alg2}.

\begin{algorithm}[t]
	\caption{Solving Support Points by CCP} 
	\label{alg2}
	\KwIn{$Z\in\mathbb{R}^{p\times n}$. }
	\textbf{Initialize:}  Sample $D^{0} =\{\mathbf{d}_t\}_{t=1}^m$ i.i.d. from $Z$, and set $j\leftarrow 1$.\\
	\While{no converged}{
		\For{$t=1,\ldots, m$ \textbf{do parallel:}}{
			Compute $\mathbf{d}_t$ by \eqref{eq:d}.\\
			Update $D^k_{j}\leftarrow \mathbf{d}_t$, and set $j \leftarrow j+1$.}
	}
\end{algorithm}

The next theorem shows that if the sequence generated by Algorithm \ref{alg2} has a limit point, then it must be a stationary point of \eqref{sp}. The following theorem can follow directly from the results of \cite[Theorem 7 ]{mak2018support}, thus we omit the proof here.
\begin{theorem}\label{th_sp}
	Let $\mathcal{D}$ be compact and convex, and $\{D^k\}$ be the sequence generated by Algorithm \ref{alg2}. If $D_0\subseteq \mathcal{D}$, any accumulation point of $\{D^k\}$ is a local minimizer of \eqref{sp}. 
\end{theorem}

Moving forward, we analyze the computational complexity of Algorithm \ref{alg2}. The primary computational cost arises from solving \eqref{eq:d}, which has a complexity of $\mathcal{O}((m+n)p)$. Thus, the overall complexity of obtaining $m$ support points is $\mathcal{O}(m(m+n)p)$.

\subsection{Online Algorithm with Dictionary Update}

In this part, we present two online algorithms, OENSC and OENSC-S, by integrating the proposed online $\ell_0$ elastic net subspace clustering method with the support point-based dictionary updating strategy.

As shown in the previous subsection, a small number of support points can adequately reflect the data distribution. Based on this insight, a straightforward idea is to employ the support points as the initial dictionary in Algorithm \ref{alg1}, as outlined in Algorithm \ref{alg4}. It is important to note that Algorithm \ref{alg4} only constructs the dictionary at the beginning by support points. Therefore, the quality of the dictionary depends on the quality of the initial data. Although this approach can reflect the data distribution to some extent, as mentioned earlier, the distribution of online data is dynamic and constantly evolving. This may limit the performance of Algorithm \ref{alg4}.

\begin{algorithm}[!t]
	\caption{OENSC Algorithm with Support Points} 
	\label{alg4}
	\KwIn{ $Z\in\mathbb{R}^{p\times n}$.}
	\textbf{Initialize:} $X^0, Y^0, U^0\in \mathbb{R}^{m\times n}, \lambda_1, \lambda_1, \sigma, \delta$, and find $ D\in \mathbb{R}^{p\times m}$ by Algorithm \ref{alg2}.\\
	
	\For{$i=1,2,\ldots $}{
		Update $(\mathbf{y}_i, \mathbf{x}_i, \mathbf{u}_i)$ by Algorithm \ref{alg1}.
	}
\end{algorithm}

To further enhance the performance of Algorithm \ref{alg4}, we present an online subspace clustering algorithm with a selective dictionary update strategy, as outlined in Algorithm \ref{alg3}. In Algorithm \ref{alg3}, lines 3 to 10 describe the dictionary update step. Notably, when the minimum distance between a new point and the existing data points exceeds a threshold $\delta$, the new point is considered an ``outlier". Ideally, we would update the dictionary for every outlier. However, the computational cost of performing frequent updates is quite high. Moreover, the effect of a small number of outlier points on the overall clustering performance is negligible, as the dictionary is constructed by support points, which effectively capture the data distribution while also being robust to outliers. To reduce computational overhead, we choose to update the entire dictionary $D$ only when the number of detected outlier points exceeds a specified threshold $mm$.

\begin{algorithm}[!t]
	\caption{OENSC Algorithm with Dictionary Update (OENSC-S)} 
	\label{alg3}
	\KwIn{ $Z\in\mathbb{R}^{p\times n}$.}
	\textbf{Initialize:} $X^0, Y^0, U^0\in \mathbb{R}^{m\times n}, D\in \mathbb{R}^{p\times m}, \lambda_1, \lambda_1, \sigma, \delta$, and set $count = 0$, $\hat{Z}  = D$.\\
	\For{$i=1,2,\ldots $}{
		\If{$\min_l \lVert \mathbf{z}_i -\mathbf{d}_l\rVert_2^2 \geq \delta^2$}{
			Set $count = count+1$;\\
			Update $\hat{Z} = [\hat{Z}, \mathbf{z}_i]$;\\
			\If{$count ~\%~ mm ==0$}{
				Update $D$ by Algorithm \ref{alg2} with $\hat{Z}$;\\
					%
				Set $count = 0$ and $\hat{Z}  = D$.
		}}	
		
		Update $(\mathbf{y}_i, \mathbf{x}_i, \mathbf{u}_i)$ by Algorithm \ref{alg1}.
	}
\end{algorithm}

In the following, we analyze the convergence of Algorithm \ref{alg4} and Algorithm \ref{alg3}. For Algorithm \ref{alg4}, it only differs from Algorithm \ref{alg1} in the dictionary initialization policy. Thus it inherits the same convergence property as Algorithm \ref{alg1}, i.e., any accumulation point of the iterates is a uniquely local minimizer of \eqref{OENSC1}. On the other hand, Algorithm \ref{alg3} introduces a selective dictionary update strategy to Algorithm \ref{alg1}. Although this modification results in a varying dictionary for each sample, the data in our algorithm is processed in a one-by-one input form. Therefore, despite the dictionary updates, any accumulation point of iterates for Algorithm \ref{alg3} still converges to a uniquely local minimizer for the vector version problem \eqref{OENSC_ve}.

Based on the previous analysis, the computational complexity of Algorithm \ref{alg4} and Algorithm \ref{alg3} are $\mathcal{O}(k_1m(p^2 + m^2))$ and $\mathcal{O}(\max\{k_1m(p^2 + m^2), k_2m(m+n)p\})$ respectively, where $m, n$, and $p$ are the number of support points, samples, and features respectively. $k_1$ and $k_2$ are iterations for Algorithm \ref{alg1} and Algorithm \ref{alg2}, respectively.

\subsection{Comparison with Online Dictionary Learning Schemes}

Over the past decades, numerous online dictionary learning methods have been developed and applied to clustering \cite{mairal2009online,dumitrescu2018dictionary,wu2021online}. For online data, these methods usually employ a two-stage updating strategy. The first stage computes the sparse representation coefficients of the data, and the second stage updates the dictionary based on these coefficients. This strategy relies on the historical information of the entire dataset, which leads to inflexible performance in coping with real-time data changes, thus making it less effective in capturing the dynamic nature of online data. Furthermore, both stages involve iterative computations, resulting in high computational complexity, which fails to meet the demand for real-time updating, and encounters significant computational bottlenecks, especially for high-dimensional data. Additionally, existing methods usually initialize the dictionary with random variables, which prevents it from effectively capturing the data distribution, thereby affecting the clustering quality.

OENSC employs a dictionary update strategy based on support points, which selects the dictionary elements according to data representativeness. Support points are chosen by minimizing the energy distance, which allows for a more accurate summary of the data distribution. This dictionary update strategy effectively adapts to changes in online data streams, adjusting the dictionary structure to reflect these variations. Furthermore, the proposed method adopts a periodic updating strategy, which makes it more flexible in coping with the drift. It effectively balances the computational overhead with the updating frequency, and is more suitable for large-scale and high-dimensional data.

\section{Experiments}\label{sec_exp}
This section conducts extensive experiments to demonstrate the effectiveness of our proposed OENSC and OENSC-S methods on six publicly available datasets, i.e., Extended Yale B (EYaleB)  \cite{lee2005acquiring}, AR \cite{martinez1998ar}, USPS \cite{hull1994database}, CIFAR-10 \cite{krizhevsky2009learning}, Network Intrusion \cite{cieslak2006combating}, and MNIST \cite{lecun1998gradient}. 
The statistics descriptions of these datasets are summarized in Table 1. 
All numerical experiments are conducted on MATLAB (R2022a) under the Windows environment with an Intel(R) Core(TM) i7 CPU @2.10 GHz and 32 GB of RAM. We begin by discussing the implementation details, followed by an analysis of the experimental results.

\begin{table}[t]
	\centering
	\caption{The statistics of the datasets.}
	\label{tb_datasets}
	\begin{tabular}{c c c c}
		\noalign{\smallskip}\hline\noalign{\smallskip}
		Datasets & Smaples 	& Dimension 		& Classes\\
		\noalign{\smallskip}\hline\noalign{\smallskip}
		EYaleB	& 2414		& 256 				& 38\\
		AR		& 2600		& 192 				& 100\\
		USPS	& 9298 		& 256				& 10\\
		CIFAR-10 & 20000	& 256 				& 10\\
		Network & 25192		& 41				& 2\\
		MNIST	& 70000		& 784 				& 10\\
		\noalign{\smallskip}\hline\noalign{\smallskip}
	\end{tabular}
\end{table}

\subsection{Experimental Settings}

\subsubsection{Compared Methods and Evaluation Criteria}
To demonstrate the competitiveness of the proposed methods, we compare them with several state-of-the-art approaches, including  SSC \cite{Elhamifar2013}, LRR \cite{Liu2013lrr}, ORSC \cite{lee2016online}, OLRC \cite{shen2016online}, OSRC \cite{chen2024online}, and TSSRC \cite{chen2023two}. Among them, SSC and LRR are batch methods, while the other four methods are designed for online data processing. Notably, OSRC and TSSRC adopt a batch-to-batch approach. The implementations of the aforementioned algorithms are publicly available in corresponding papers, and we use them directly in our experiments without revision. For a fair comparison, the codes used in our experiments are provided by the authors' websites or other publicly accessible sources, and the involved parameters of these competitors are carefully tuned according to the respective literature. In OENSC, the dictionary is only updated initially, whereas in OENSC-S, the dictionary is continuously updated using data flow, with both methods utilizing support points for updating the dictionary. 

To quantify the clustering performance, three widely used clustering performance metrics \cite{christopher2008introduction}, i.e., clustering accuracy (ACC), normalized mutual information (NMI), and purity, are reported in our experiments. The higher the value for these criteria, the better the performance. 

\begin{table*}[!t]
	\centering
	\caption{The clustering performance (\%) comparison.}
	\label{tb_exp}
	\small
	\begin{tabular}{l l c c c c c c c c }
		\noalign{\smallskip}\hline\noalign{\smallskip}
		\multirow{2}{*}{Datasets}	& \multirow{2}{*}{Metrics} &  \multicolumn{8}{c}{Methods} \\
		\cmidrule(lr){3-10}
		& 	& SSC		& LRR		& ORSC		& OLRSC		& OSRC	& TSSRC	& OENSC & OENSC-S \\
		
		\noalign{\smallskip}\hline\noalign{\smallskip}
		\multirow{2}{*}{EYaleB} & ACC	 & 61.18 &  65.37  & 68.72 & 70.26
		& 73.55 & 69.10 & \underline{75.38} & \textbf{76.97}  \\
		& NMI   & 67.13 & 68.75 & 70.85 & 71.73 & {69.65} & 73.90 & \underline{75.56} & \textbf{77.38}\\
		& Purity   & 62.92 & 68.95 & 64.93 & 66.58 & {71.67} & 71.96 & \underline{72.42} & \textbf{76.06}\\
		
		\noalign{\smallskip}\hline\noalign{\smallskip}
		\multirow{2}{*}{AR} & ACC   & 64.75 & 63.16 & 67.51 & 73.15 & 72.71 & 74.10 & \underline{74.23} & \textbf{75.08}\\
		& NMI   & 67.13 & 68.75 & 70.84 & \textbf{78.46} & {68.25} & 69.95 & {73.47} & \underline{74.42}\\
		& Purity   & 65.28 & 66.58 & 69.42 & \textbf{74.83} & {69.27} & 67.56 & {71.79} & \underline{73.63}\\
		
		\noalign{\smallskip}\hline\noalign{\smallskip}
		\multirow{2}{*}{USPS} 	& ACC   & 59.14 & 57.79  & 67.51 & 67.46 & 70.24 & 71.73 & \underline{72.98} & \textbf{73.22}\\
		& NMI   & 62.96 & 63.54 & 65.56 & 66.70 & 69.33 & 68.57 & \textbf{75.60} & \underline{75.05}\\
		& Purity   & 61.27 & 61.04 & 61.03 & 66.73 & 64.67 & 61.31 & \underline{66.91} & \textbf{68.14}\\
		
		\noalign{\smallskip}\hline\noalign{\smallskip}
		\multirow{2}{*}{CIFAR-10} & ACC   & - & - & 45.36 & 44.90 & 56.69 & 56.51 & \textbf{57.89}& \underline{57.35}\\
		& NMI   & - & - & 47.16 & 52.49 & 49.86 & 47.69 & \textbf{55.03} & \underline{54.60} \\
		& Purity   & - & - & 41.67 & 41.84 & 44.82 & 46.71 & \underline{48.62} & \textbf{51.19} \\
		
		\noalign{\smallskip}\hline\noalign{\smallskip}
		\multirow{2}{*}{Network} & ACC   & - & - & 76.74 & 80.99 & 86.69 & 84.95 & \underline{87.76} & \textbf{88.54}\\
		& NMI   & - & - & 77.43 & 83.97 & 85.64 & 84.26 & \underline{88.06} & \textbf{91.84}\\
		& Purity   & - & - & 79.81 & 85.15 & {89.32} & {88.37} & \underline{90.27} & \textbf{92.69}\\

		\noalign{\smallskip}\hline\noalign{\smallskip}
		\multirow{2}{*}{MNIST} & ACC    & - & - & 48.09 & 49.14 & 50.40 & {52.45} & \underline{57.39} & \textbf{58.81}\\
		& NMI   & - & - & 41.18 & 44.83 & 48.16 & 47.35 & \underline{52.19} & \textbf{54.27} \\
		& Purity   & - & - & 45.28 & 45.12 & 54.63 & 54.82 & \textbf{60.91} & \underline{59.74}\\
		
		\noalign{\smallskip}\hline\noalign{\smallskip}
	\end{tabular}
\end{table*}

\subsubsection{Algorithm Settings}
The proposed $\ell_0$-ENSC model \eqref{OENSC_ve} involves two tuning parameters, i.e., $\lambda_1$ and $\lambda_2$, which should be carefully selected through the grid search strategy over the candidate set $\lambda_1 \in \{10^{-5}, 5\times 10^{-5}, \ldots, 1\}$ and $\lambda_2 \in \{2^{-7}, 2^{-6}, \ldots, 2^7\}$. By fixing each parameter in turn, the numerical performance of each parameter combination can be evaluated and utilized as a basis for parameter selection. In the numerical experiments, the number of support points $m$ is set to $10c$, where $c$ is the number of categories, except for AR ($m=200$) and Yale ($m=5c$). The dictionary update period $mm$ is set to 100. 
For the proposed OENSC and OENSC-S, the iterative process is terminated when the generated sequence
$\{\mathbf{x}_i^k, \mathbf{y}_i^k\}$ satisfies
\begin{eqnarray*}
	\max\left\{ \frac{ \|\mathbf{x}_i^{k+1} - \mathbf{y}_i^{k+1}\|_2}{1 + \|\mathbf{y}^{k+1}\|_2},  \frac{\| \mathbf{x}_i^{k+1} - \mathbf{x}_i^{k}\|_2}{1 + \|\mathbf{x}_i^k\|_2} \right\} \leq 10^{-3},
\end{eqnarray*}
or when the number of iterations reaches 100. To mitigate the impact of randomness due to initialization, each experiment is repeated 30 times.

 \begin{figure}[t]
 	\centering
 	\subfigure[USPS]{\includegraphics[width=0.45\hsize]{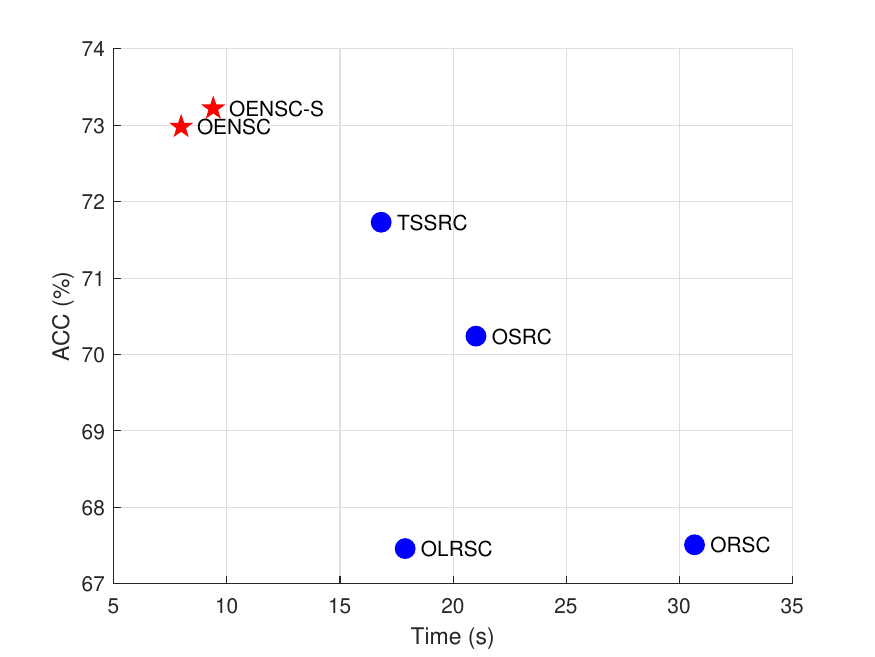}} 
 	\subfigure[CIFAR-10]{\includegraphics[width=0.45\hsize]{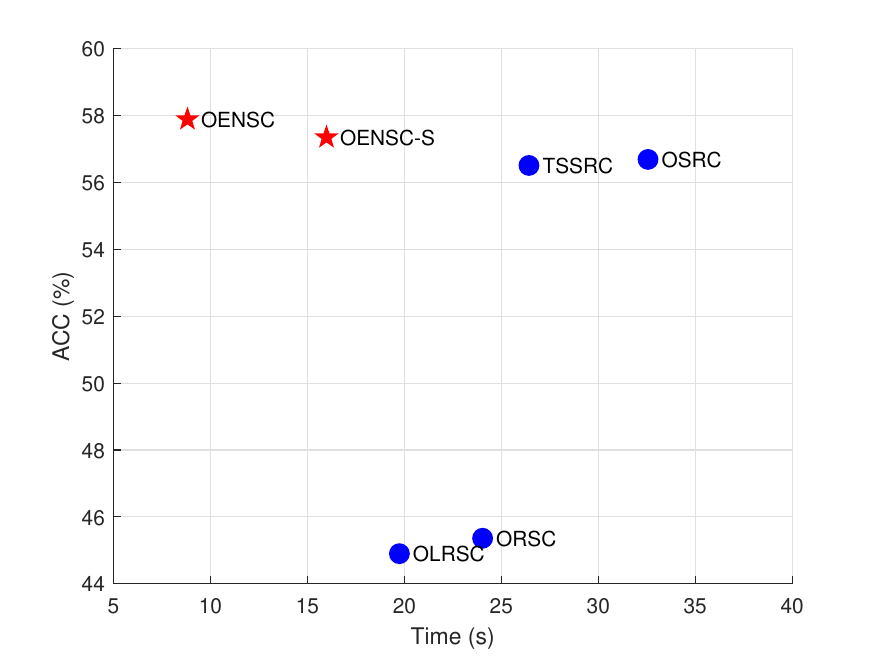}}\\
 	\subfigure[Network]{\includegraphics[width=0.45\hsize]{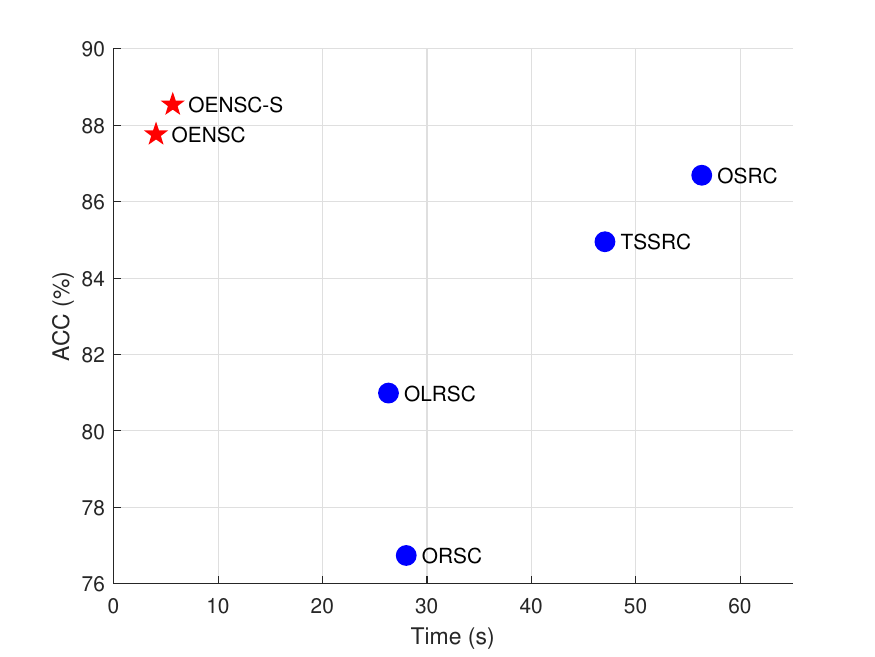}}
 	\subfigure[MNIST]{\includegraphics[width=0.45\hsize]{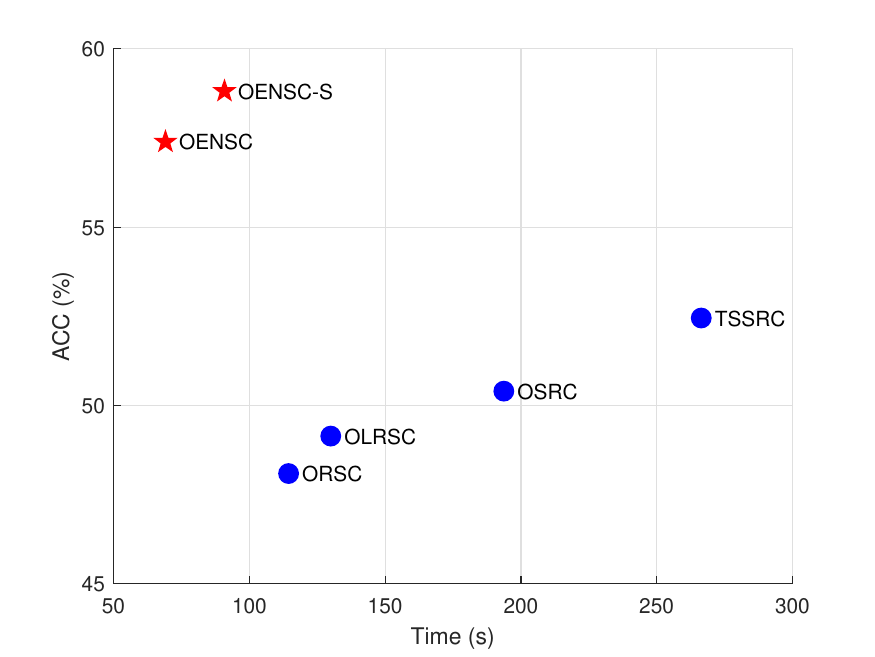}}
	\caption{The ACC versus running times on different datasets.}
 	\label{fig_acc_times}
 \end{figure}

\begin{figure}[t]
	\centering
	\subfigure[ACC on EYaleB]{\includegraphics[width=0.3\hsize]{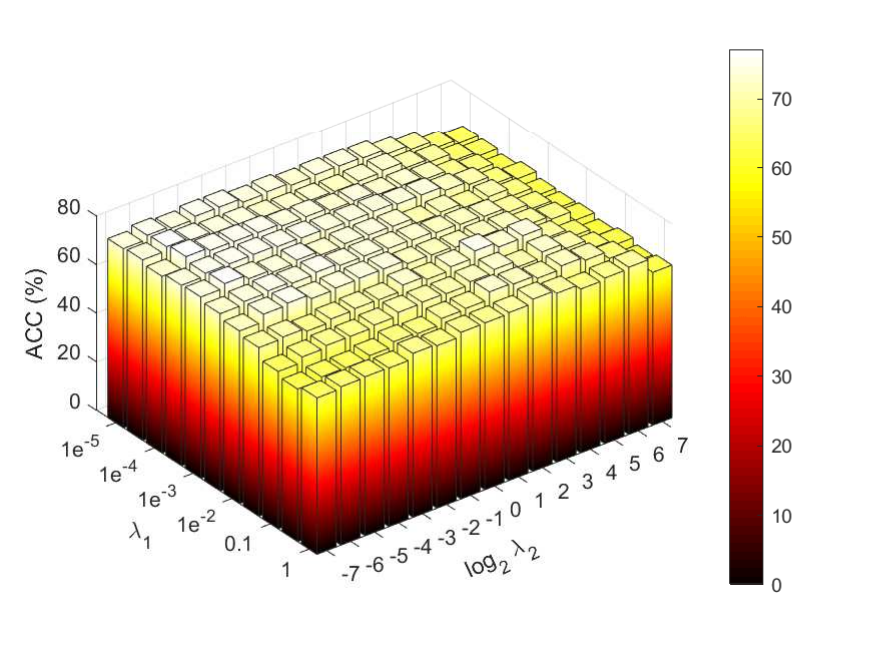}} 
	\subfigure[NMI on EYaleB]{\includegraphics[width=0.3\hsize]{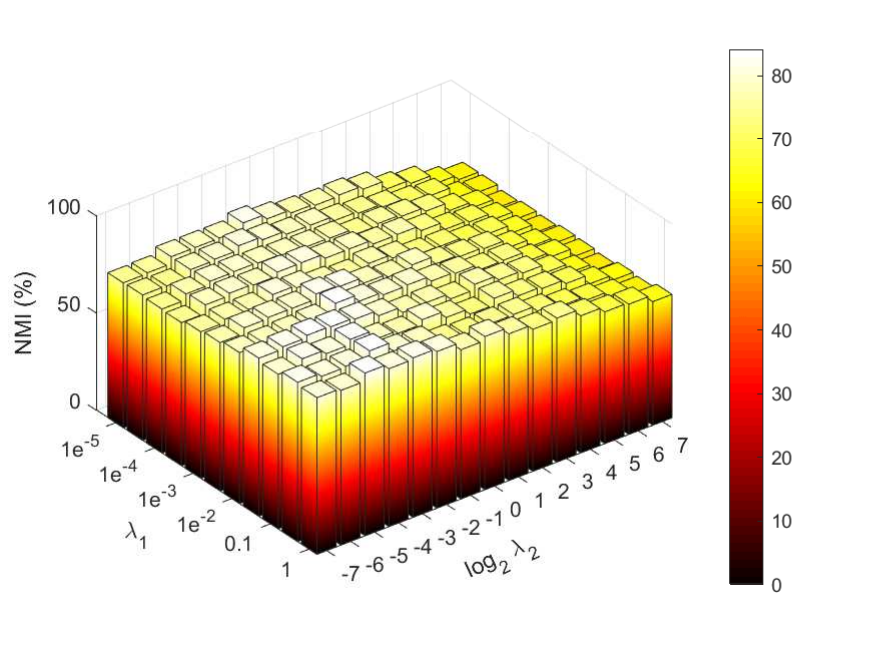}}
	\subfigure[Purity on EYaleB]{\includegraphics[width=0.3\hsize]{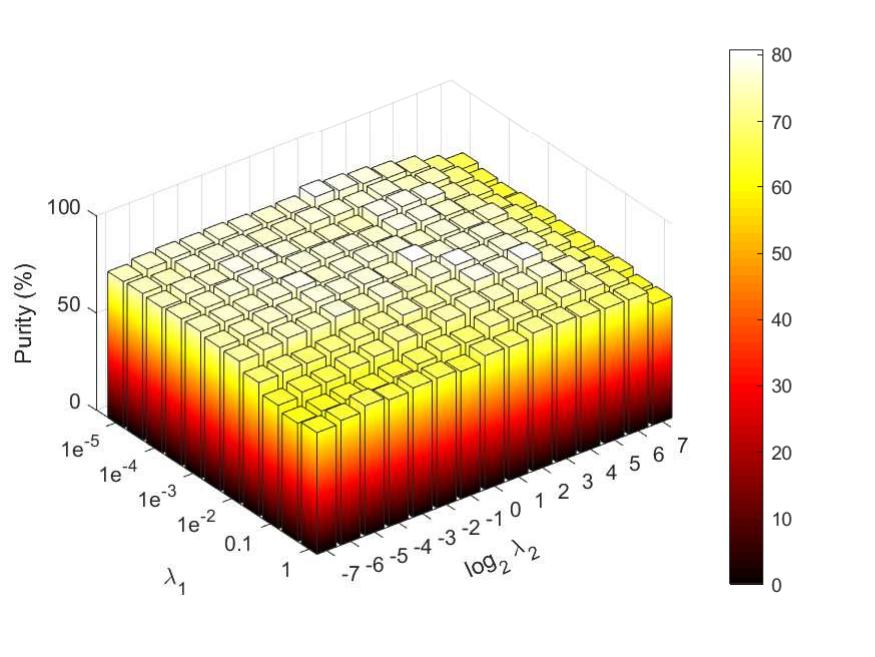}}\\
	\subfigure[ACC on USPS]{\includegraphics[width=0.3\hsize]{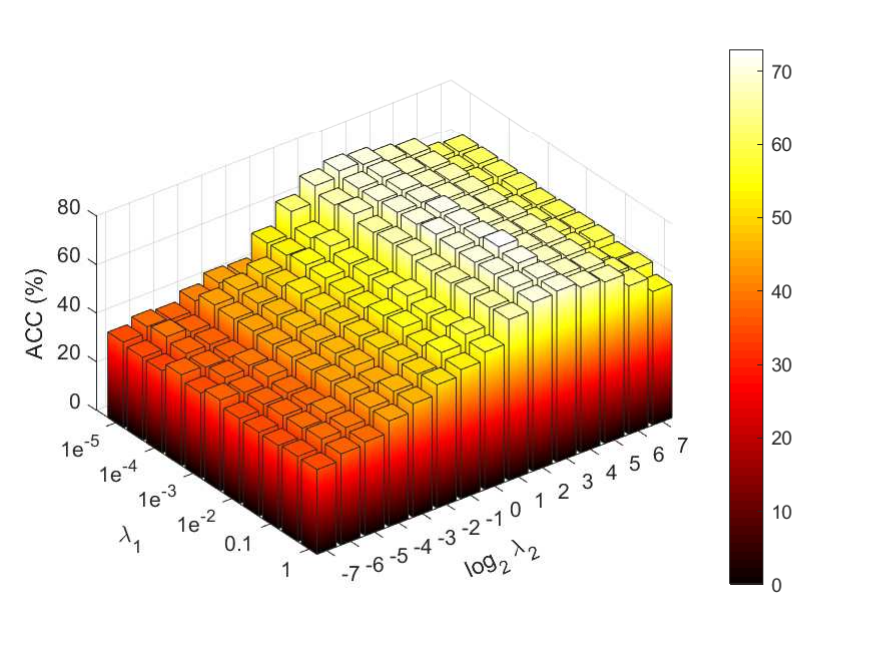}} 
	\subfigure[NMI on USPS]{\includegraphics[width=0.3\hsize]{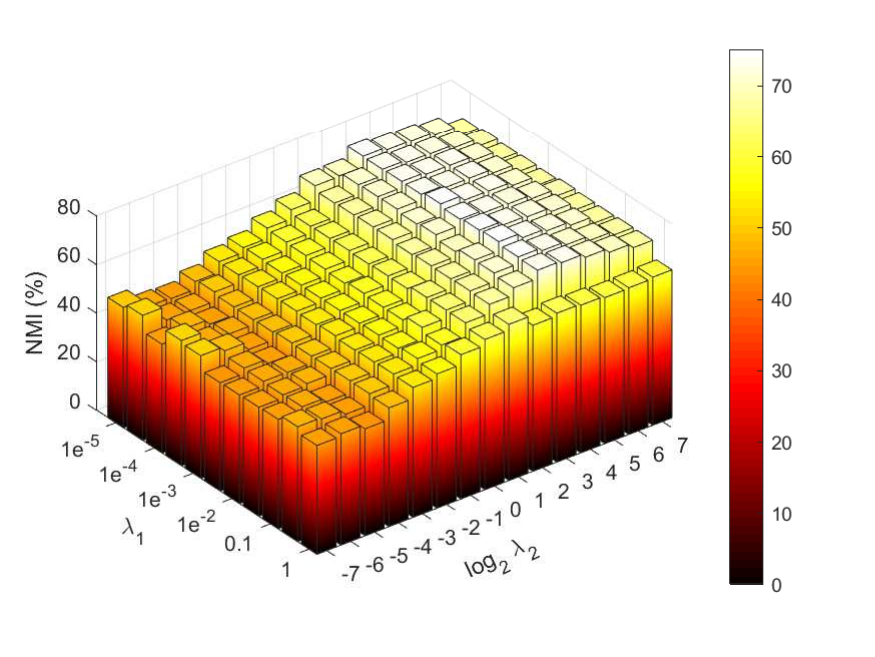}}
	\subfigure[Purity on USPS]{\includegraphics[width=0.3\hsize]{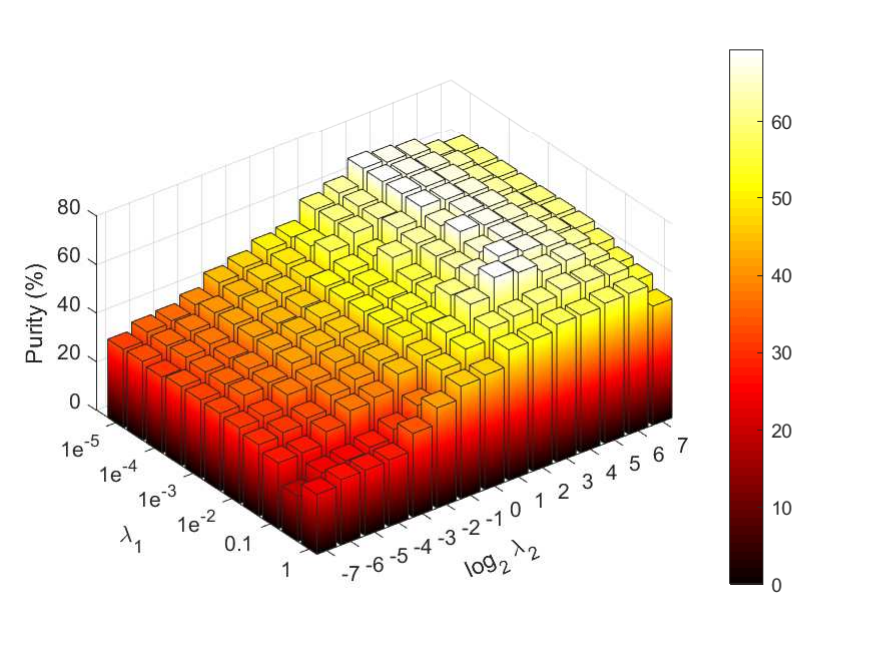}}
	\caption{The effect of tuning parameter on EYaleB and USPS.}
	\label{fig_par}
\end{figure}

\begin{figure}[t]
	\centering
	\subfigure[SSC, ACC = 67.33\%]{\includegraphics[width=0.3\hsize]{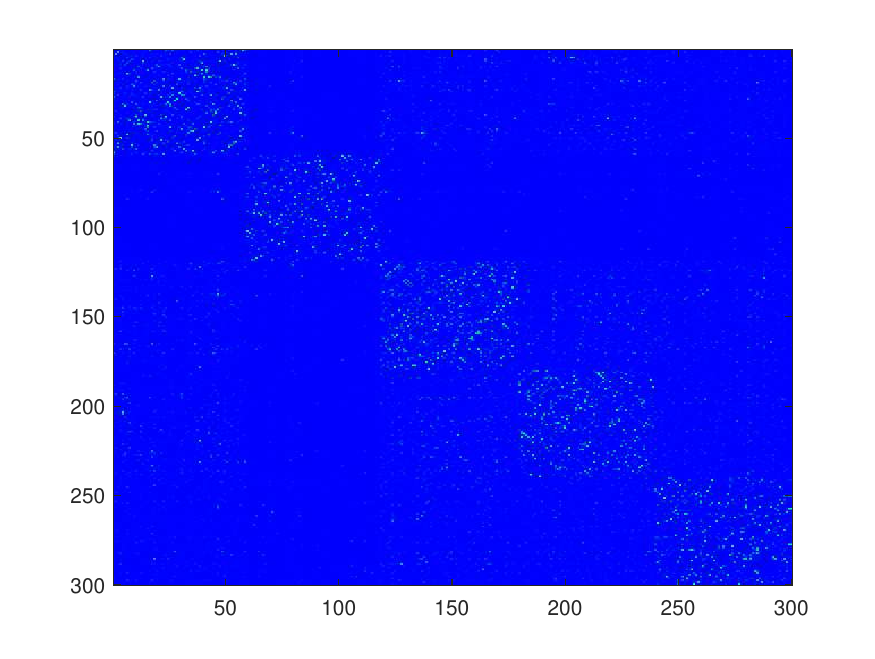}}
	\subfigure[LRR,  ACC = 62.36\%]{\includegraphics[width=0.3\hsize]{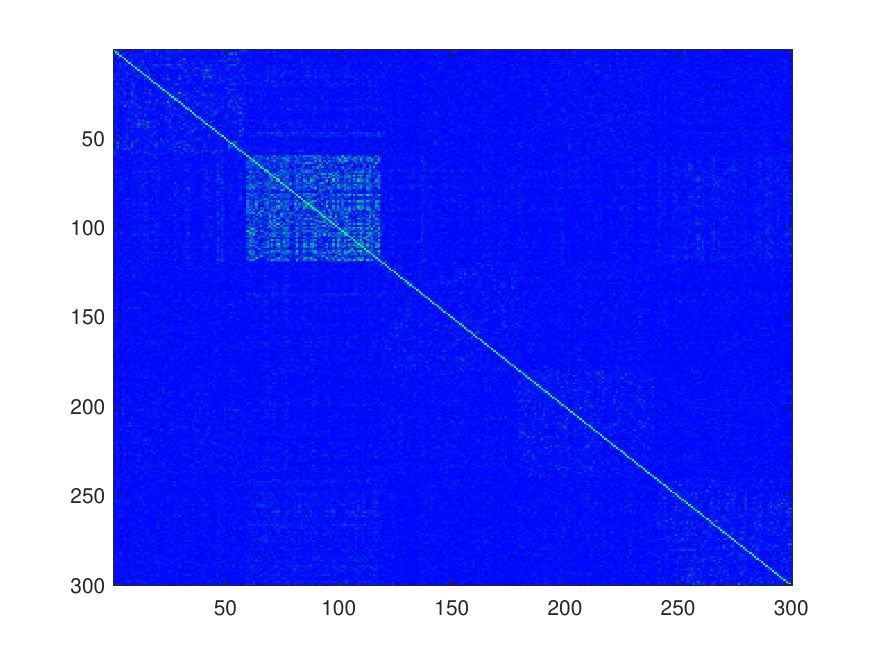}}
	\subfigure[OLRSC,  ACC = 93.58\%]{\includegraphics[width=0.3\hsize]{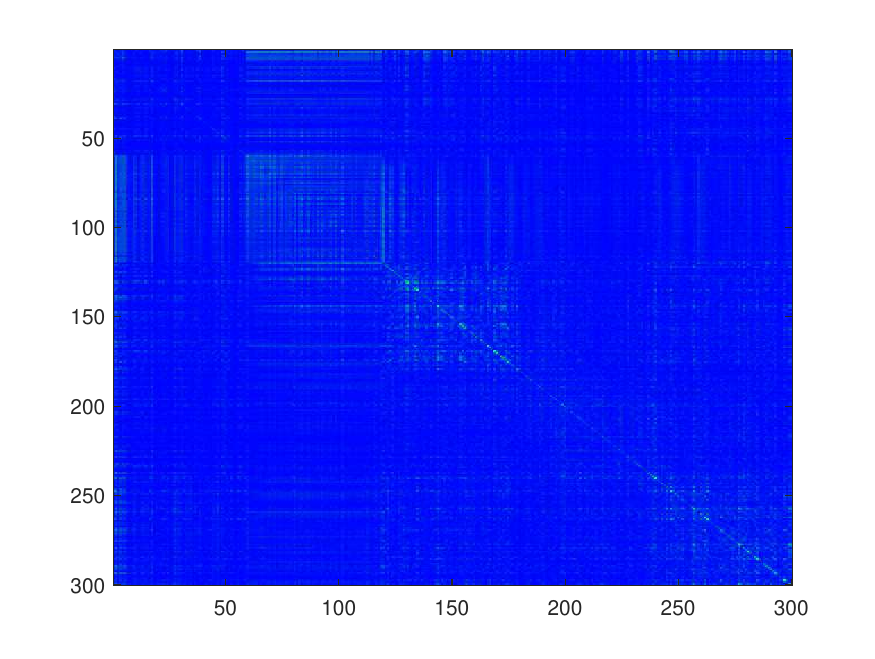}}
	\\
	\subfigure[OSRC,  ACC = 88.33\%]{\includegraphics[width=0.3\hsize]{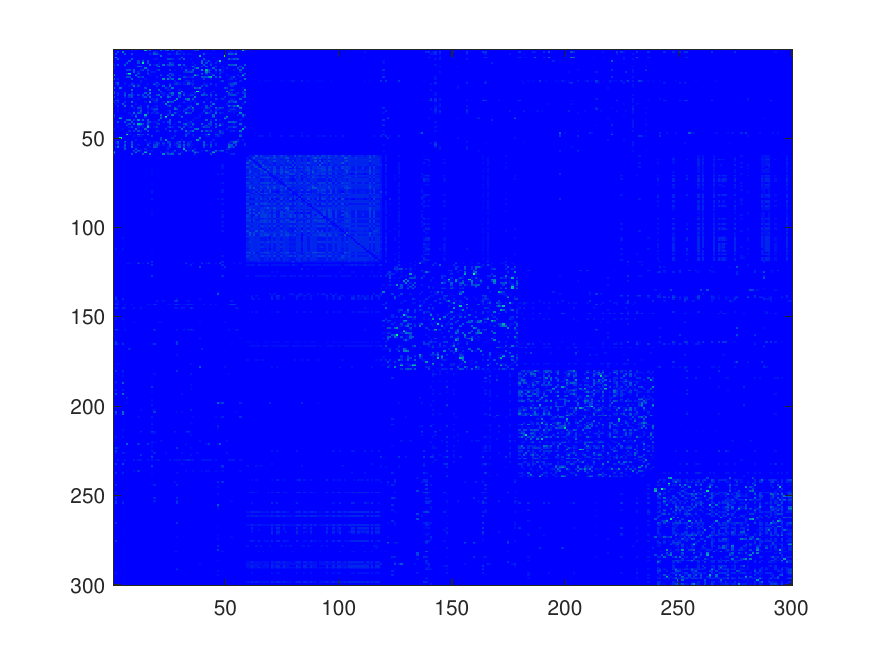}}
	\subfigure[TSSRC, ACC =  94.67\%]{\includegraphics[width=0.3\hsize]{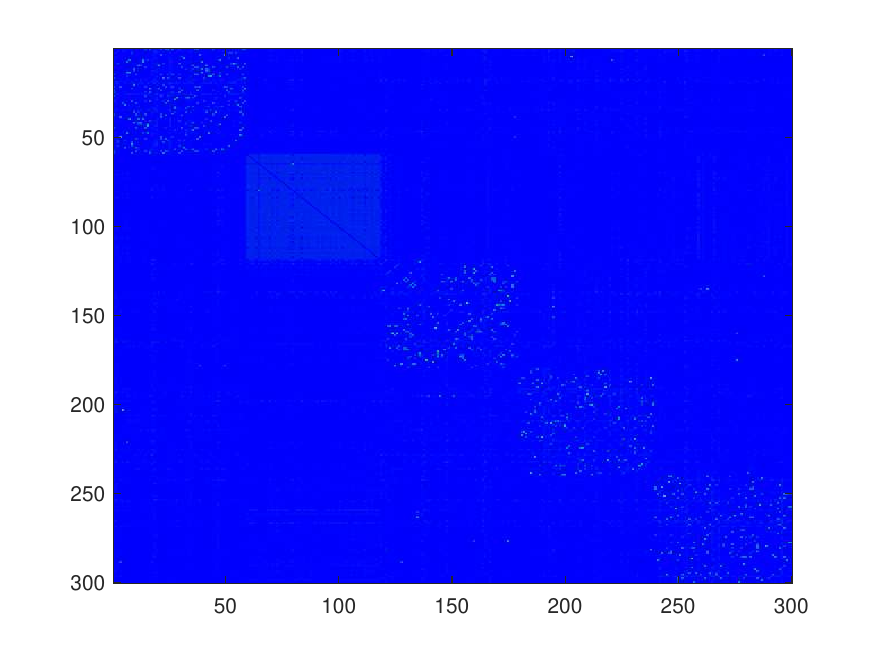}}
	\subfigure[OENSC, ACC = 96.70\%]{\includegraphics[width=0.3\hsize]{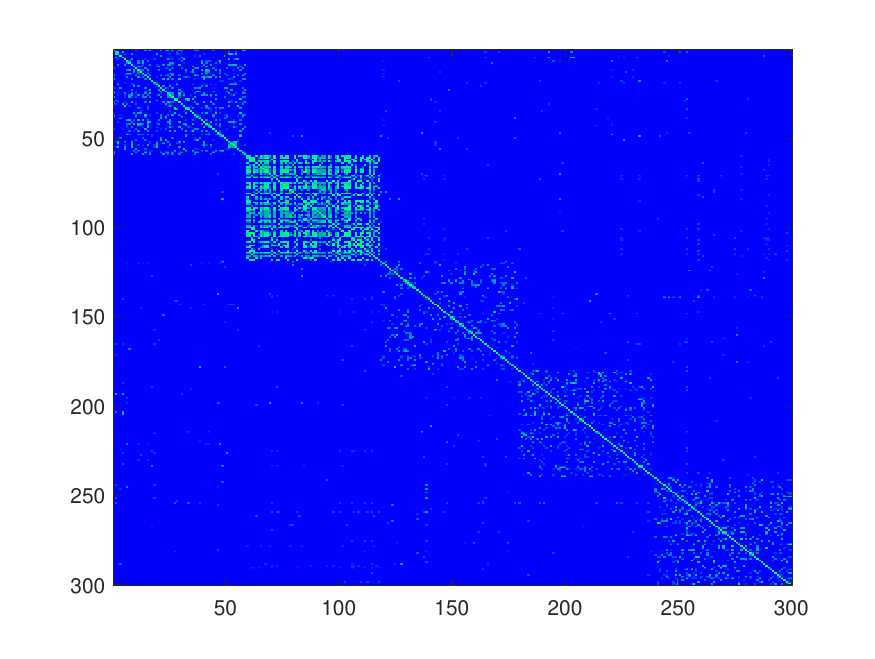}} \\
	\caption{The visualization of coefficient matrices on USPS with 300 samples drawn from 5 classes.}
	\label{fig_HEAT}
\end{figure}

\begin{figure}[t]
	\centering
	\subfigure[ACC (Salt)]{\includegraphics[width=0.3\hsize]{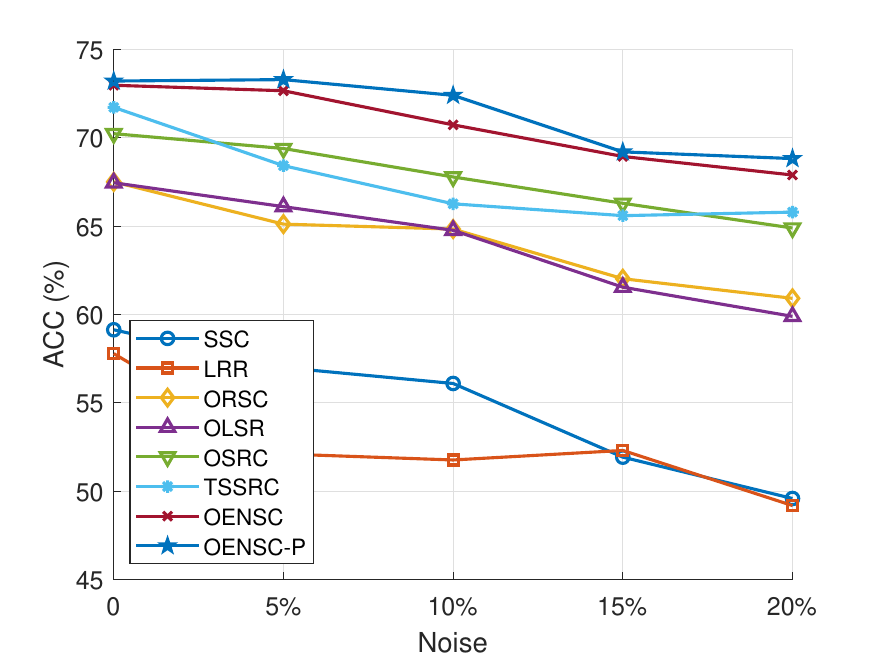}}
	\subfigure[NMI (Salt)]{\includegraphics[width=0.3\hsize]{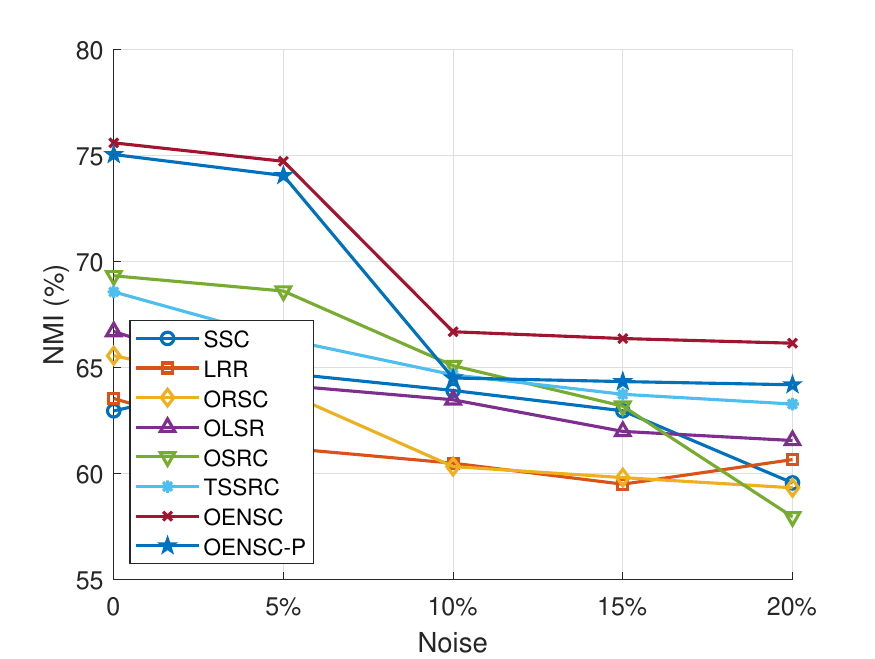}}
	\subfigure[Purity (Salt)]{\includegraphics[width=0.3\hsize]{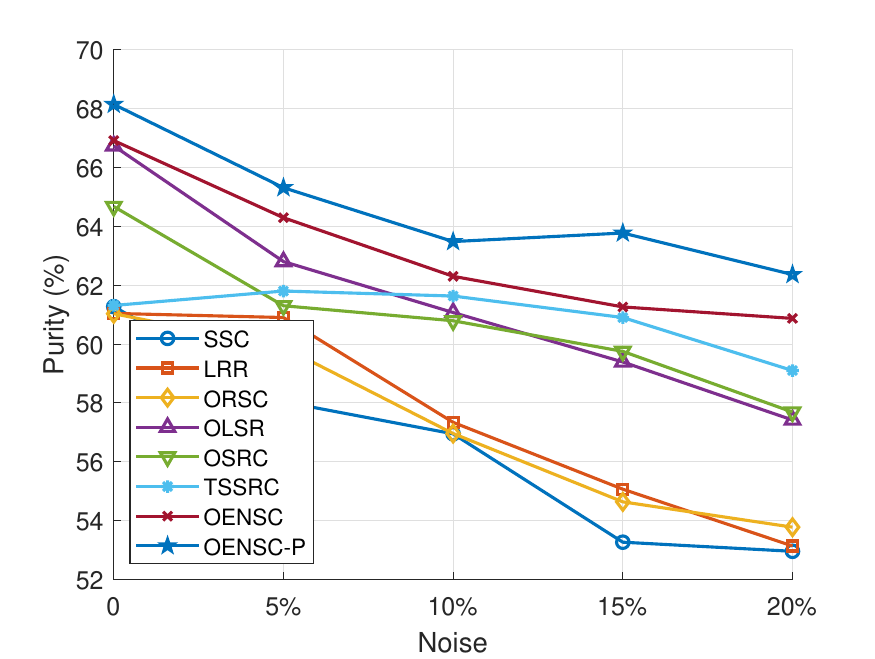}}\\
	\subfigure[ACC (Speckle)]{\includegraphics[width=0.3\hsize]{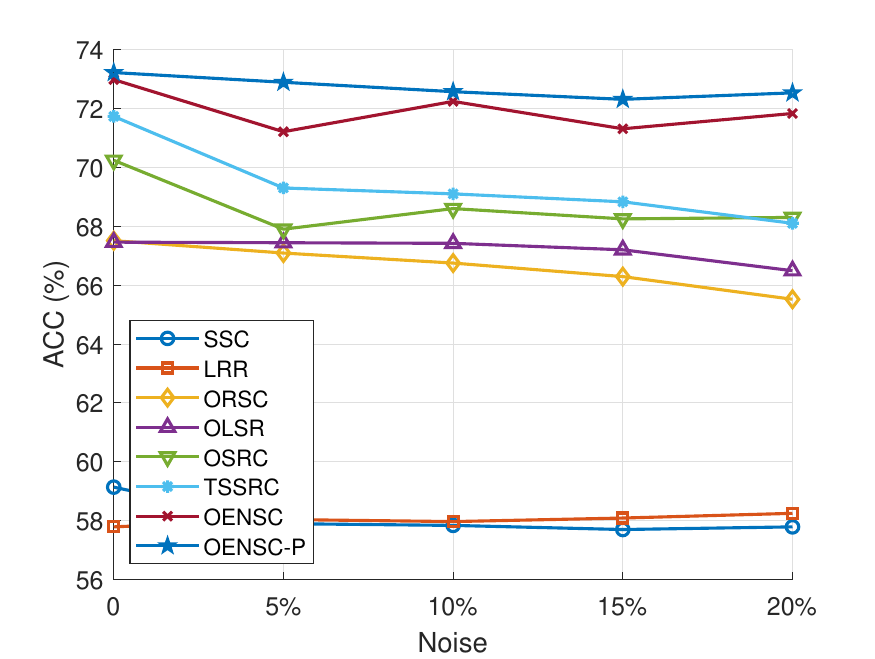}}
	\subfigure[NMI (Speckle)]{\includegraphics[width=0.3\hsize]{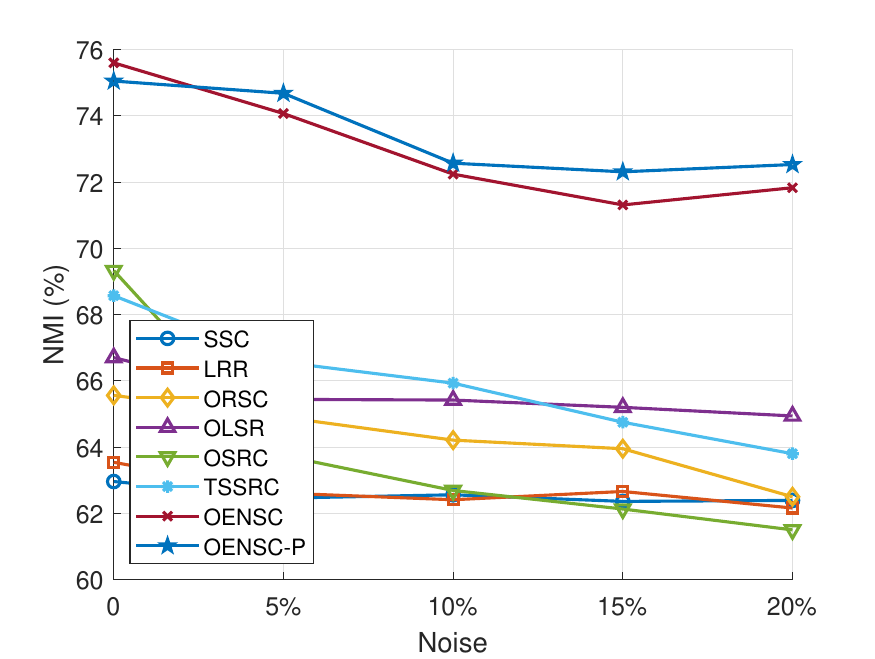}}
	\subfigure[Purity (Speckle)]{\includegraphics[width=0.3\hsize]{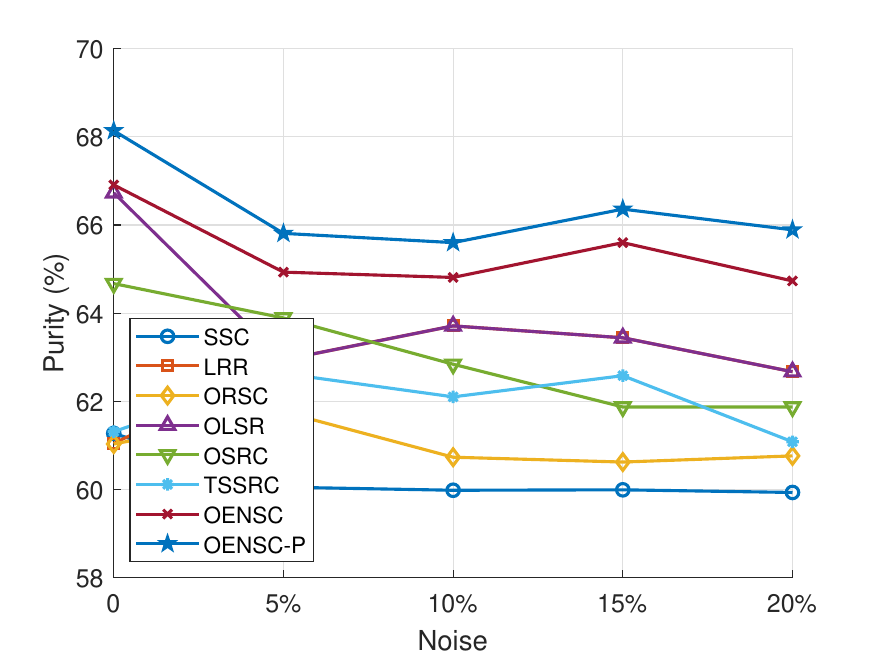}}\\
	\caption{The clustering performance on USPS with different types of noise.}
	\label{fig_noise}
\end{figure}

\begin{figure}[t]
	\centering
	\subfigure[USPS]{\includegraphics[width=0.45\hsize]{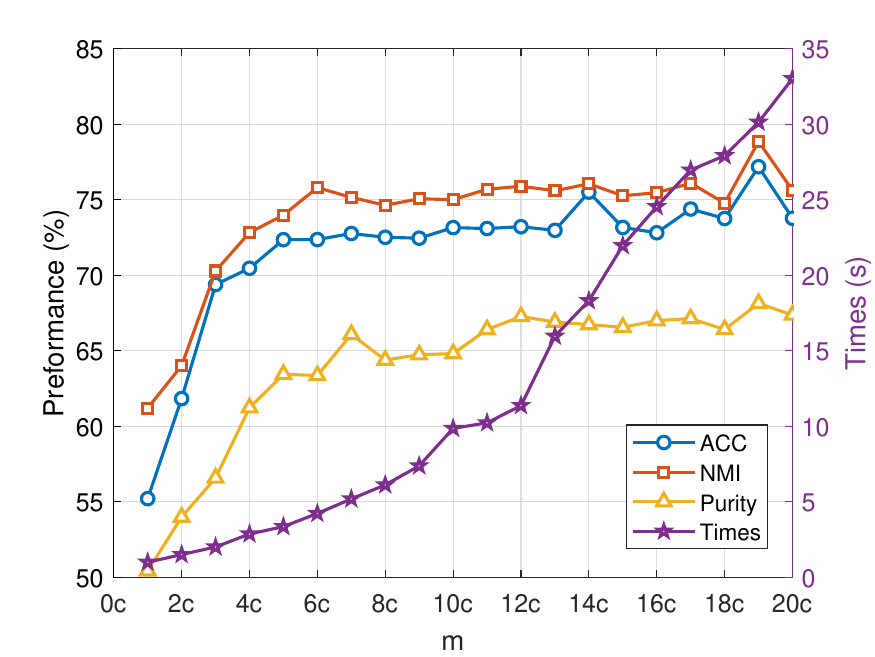}} 
	\subfigure[Network]{\includegraphics[width=0.45\hsize]{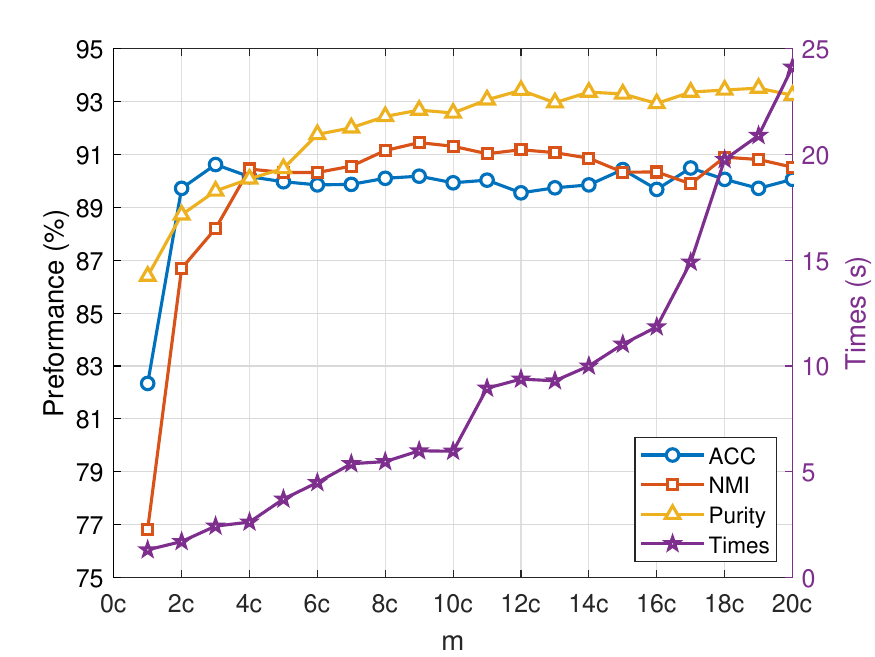}}
	\caption{The effect for the number of support points.}
	\label{fig_change_m}
\end{figure}

\subsection{Experimental Results}

Table \ref{tb_exp} summarizes the clustering performance of all compared methods,  and Figure \ref{fig_acc_times} reports the ACC versus running times. The best results are highlighted in bold, while the second-best are underlined. From these tables, it can be seen that the following observations
\begin{itemize}
	\item For all selected datasets, the proposed OENSC and OENSC-S consistently outperform other compared methods, except on the AR, where the NMI and Purity metrics are slightly lower than those of OLRSC. Specifically, OENSC-S improves the NMI by approximately 6\% on USPS and Network datasets compared to other methods.
	\item OENSC-S generally delivers better clustering performance than OENSC due to frequent dictionary updates. However, this comes at a cost of higher computation time. For instance, OENSC-S improves ACC by about 1.5\% on the EYaleB and MNIST datasets, but its runtime is much longer than that of OENSC. This highlights a trade-off between clustering performance and computational efficiency, with OENSC being more suitable for applications requiring a faster response.	
	\item Both OENSC and OENSC-S are computationally efficient compared to other methods, especially on large-scale datasets. For example, on the Network dataset, OENSC requires only 2.92 seconds, whereas methods like OSRC and TSSRC require over 60 seconds. This makes OENSC and OENSC-S more suitable for scenarios requiring both speed and accuracy.
\end{itemize}

\subsection{Discussion}
\subsubsection{Parameter Analysis}

Note that in model \eqref{OENSC_ve}, there are two tuning parameters $\lambda_1$ and $\lambda_2$. 
Figure \ref{fig_par} shows the effect of these parameters on clustering performance for the EYaleB and USPS datasets. From these visualizations, the following conclusions can be drawn
\begin{itemize}
	\item These parameters do have an impact on the numerical performance of OENSC-S, but the extent of their impact differs across datasets. Generally, $\lambda_2$ has a significant effect, as the $\ell_2$ regularization promotes to group the highly correlated samples into the same cluster. As shown in the figure, the performance metrics, particularly on the USPS, exhibit more variability as $\lambda_2$ changes.
	\item Although ACC does not exhibit a uniform and pronounced trend as  $\lambda_1$ varies, fully utilizing the sparsity of data can effectively enhance computational efficiency.
\end{itemize}

\subsubsection{Visual Analysis}

Figure \ref{fig_HEAT} shows the coefficient matrix for the sub-dataset of USPS with 60 samples in each of the 5 categories.  It is evident that the coefficient matrix of OENSC is notably sparse, which indicates that combining $\ell_0$ norm with the Frobenius norm allows for flexible selection. 
Furthermore, the coefficient matrix of OENSC exhibits a more explicit block-diagonal structure compared to other methods, aligning with Theorem \ref{th_bdq}. This highlights the ability of the proposed model to identify and select meaningful samples from the raw data for self-expression, thereby enhancing the clustering performance.

\subsubsection{Robustness Analysis}

In this part, we investigate the robustness of compared methods to different noises. To demonstrate the robustness of the proposed method, we conduct a comprehensive study on USPS with various two types of non-Gaussian noise, including
\begin{itemize}
	\item Type-I: Salt-and-pepper noise with percentages $\{ 5\%, 10\%, 15\%, 20\%\}$.
	\item Type-II: Speckle noise with percentages $\{ 5\%, 10\%, 15\%, 20\%\}$.
\end{itemize}

Figure \ref{fig_noise} presents the classification accuracy. It can be concluded that our OENSC and OENSC-S maintain a significant performance advantage over the other methods as the noise level increases. 
In particular, for Speckle noise, OENSC-S maintains more than 70\% for ACC and NMI even at 20\% noise, whereas the performance of the other methods drops sharply. These results highlight that the proposed methods not only achieve the best performance in clean data scenarios, but also maintain their advantage under significant noise, making them well-suited for real-world applications where data is often contaminated with different types of noise.

\begin{table}[t]
	\centering
	\caption{The ablation study of OENSC-S on USPS and Network.}
	\label{tb_ablation}
	\begin{tabular}{c c c c c c c c}
		\noalign{\smallskip}\hline\noalign{\smallskip}
		\multirow{2}{*}{$\|\mathbf{x}_i\|_0$} & \multirow{2}{*}{$\|\mathbf{x}_i\|_2^2$} & \multicolumn{3}{c}{USPS} & \multicolumn{3}{c}{Network}\\
		\cmidrule(lr){3-5}\cmidrule(lr){6-8}
		& 				& ACC & NMI & Purity & ACC & NMI & Purity \\
		\noalign{\smallskip}\hline\noalign{\smallskip}
		$\checkmark$ & $\checkmark$ & 72.98 & 70.91 & 76.32 & 88.30 & 90.28 & 89.05 \\
		$\checkmark$ & $\times$ 	& 53.36 & 64.73 & 54.69 & 67.31 & 73.39 & 68.98 \\
		$\times$ 	 & $\checkmark$	& 69.28 & 66.10 & 72.49 & 73.81 & 75.31 & 70.87 \\
		\noalign{\smallskip}\hline\noalign{\smallskip}
	\end{tabular}
\end{table}

\subsubsection{Support Points Analysis}

In this part, we aim to analyze the effect of the number of support points ($m$) on the clustering performance. Although the larger $m$ can more accurately capture the data distribution, it also results in increased computational costs. Figure \ref{fig_change_m} illustrates the clustering performance and computation time curves of our OENSC as $m$ varies on the USPS and Network datasets. From this figure, it can be observed that the clustering performance initially improves and then stabilizes as $m$ increases, whereas the computation time exhibits a rapid increase as $m$ grows.

\subsubsection{Ablation Study}

To gain a further understanding of the proposed method, we evaluate the effect of the proposed $\ell_0$-ENSC model via an ablation study on USPS and Network datasets. 
Specifically, we compare the numerical performance of removing the $\ell_0$ term and the $\ell_2$ term, respectively, as shown in Table \ref{tb_ablation}.  It is observed that the ACC, NMI, and Purity values of the proposed method are almost 4\%, 5\%, and 4\% higher for the USPS dataset, respectively, and these values are around 15\%, 15\%, and 18\% higher for the Network dataset, respectively. 
Additionally, it indicates that the $\ell_2$ term contributes more significantly to the performance, due to the ability to encourage group effects. 
Overall, the ablation study demonstrates that the inclusion of $\ell_0$ term and $\ell_2$ term by our method can effectively and robustly enhance the subspace clustering performance.

\begin{figure}[t]
	\centering
	\subfigure[EYaleB]{\includegraphics[width=0.3\hsize]{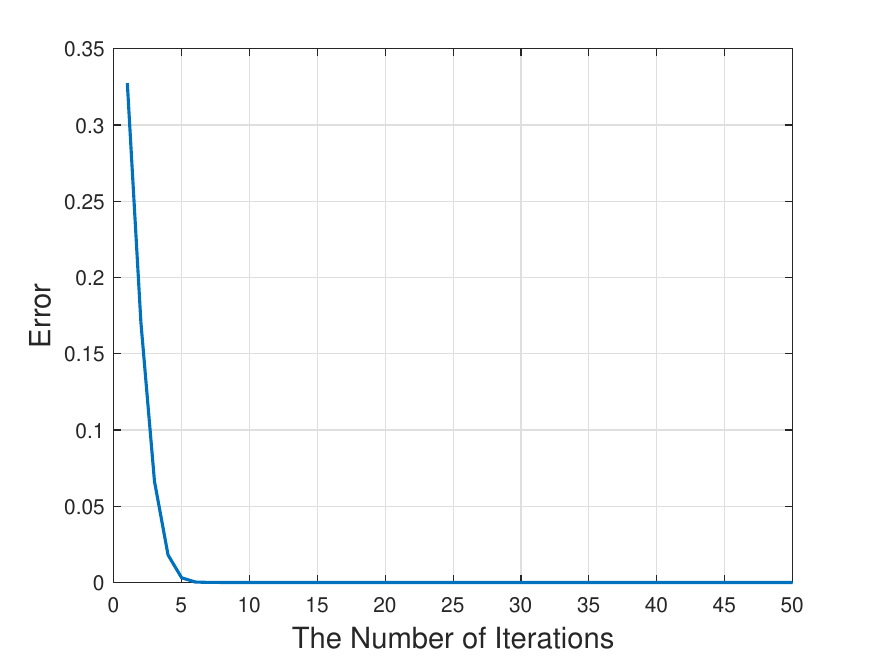}} 
	\subfigure[AR]{\includegraphics[width=0.3\hsize]{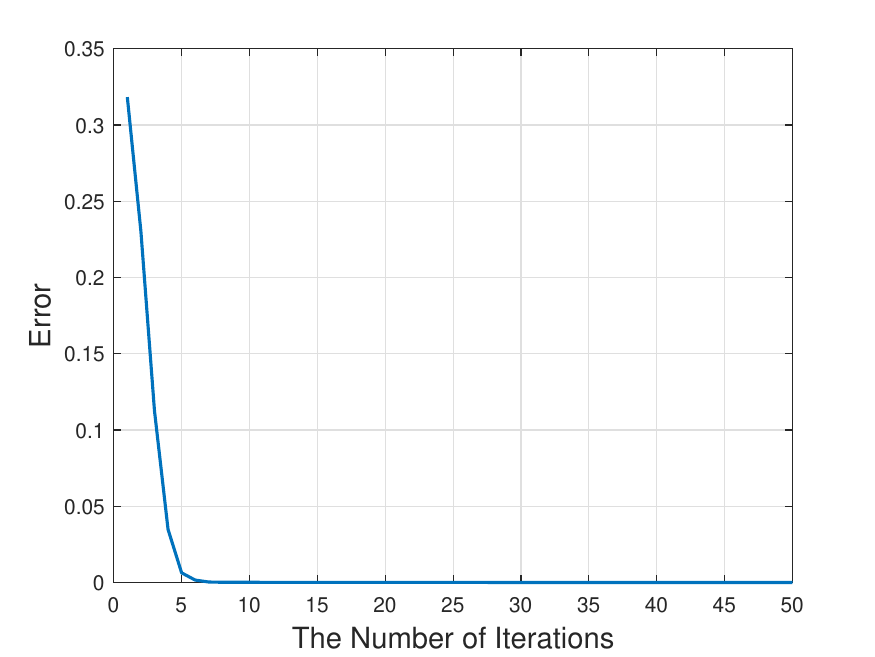}}
	\subfigure[USPS]{\includegraphics[width=0.3\hsize]{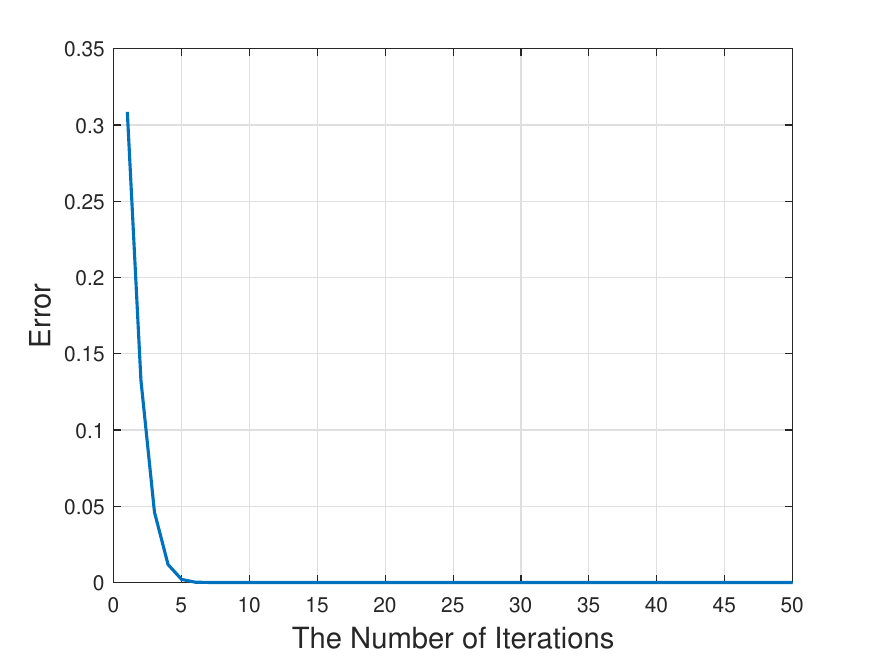}}\\
	\subfigure[CIFAR-10]{\includegraphics[width=0.3\hsize]{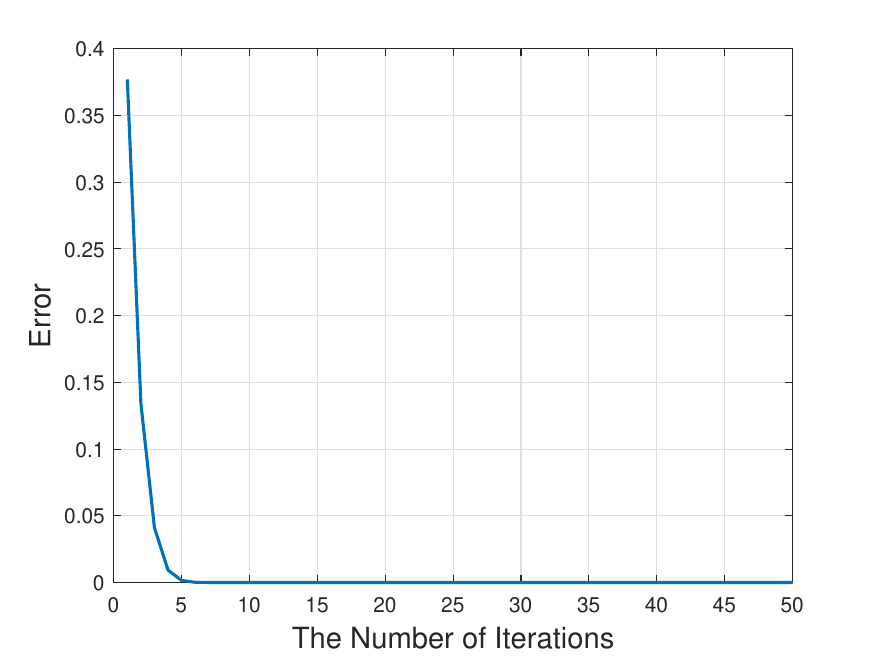}}
	\subfigure[Network]{\includegraphics[width=0.3\hsize]{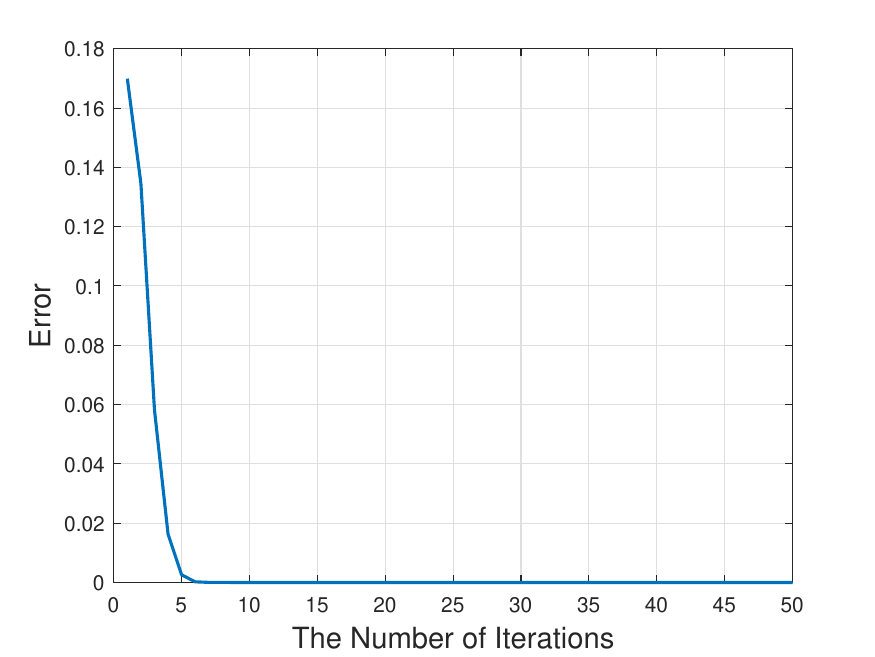}}
	\subfigure[MNIST]{\includegraphics[width=0.3\hsize]{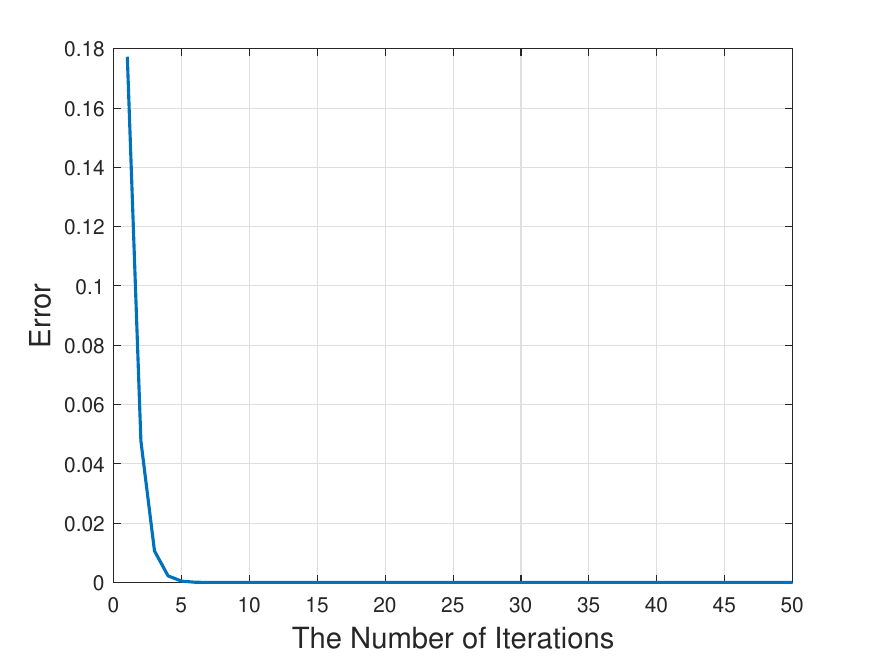}}
	\caption{The error versus the number of iterations.}
	\label{fig_error}
\end{figure}

\begin{table*}[t]
	\centering
	\caption{The runtime(s) of all compared methods.}
	\label{tb_times}
	\small
	\begin{tabular}{l c c c c c c c c }
		\noalign{\smallskip}\hline\noalign{\smallskip}
		\multirow{2}{*}{Datasets}	& \multicolumn{8}{c}{Methods} \\
		\cmidrule(lr){2-9}
		& SSC		& LRR		& ORSC		& OLRSC		& OSRC	& TSSRC	& OENSC & OENSC-S \\
		
		\noalign{\smallskip}\hline\noalign{\smallskip}
		{EYaleB} & 61.59	& 33.34 & 11.22  & 10.05  & 2.58  & \underline{2.16} & \textbf{1.23} & 2.54 \\ 
		
		{AR} & 125.20 & 31.71 & 12.83  & 12.46  & 8.93  & \underline{8.59} & \textbf{7.42} & 11.02 \\ 
		
		{USPS} 	& 1841.10 & 140.32 & 30.67  & 17.88  & 21.01  & 16.82 & \textbf{7.98} & \underline{9.40} \\ 
		
		{CIFAR-10} & - & - & 24.02 & 19.73  & 32.55 & 26.41 & \textbf{8.80} & \underline{15.97} \\ 
		
		{Network} & - & - & 28.01 & 26.31  & 56.31 & 47.04 & \textbf{4.05} & \underline{5.65} \\ 
		
		{MNIST} & - & - & 114.40 & 129.93 & 193.68  & 266.39 & \textbf{69.05} & \underline{90.78} \\ 
		
		\noalign{\smallskip}\hline\noalign{\smallskip}
	\end{tabular}
\end{table*}

\subsubsection{Convergence Analysis}

In Section \ref{sec_alg}, we have theoretically demonstrated that any accumulation point of the sequence generated by Algorithm \ref{alg3} is a P-stationary point of the proposed model \eqref{OENSC_ve}. To illustrate the convergence behavior, Figure \ref{fig_error} shows the error variation curve.  Clearly, the error decreases rapidly in the first few iterations and stabilizes within 10 iterations, which guarantees the efficiency of OENSC-S in practical applications.

Table \ref{tb_times} provides the runtime in seconds of all compared methods. It can be concluded that the computational efficiency of the proposed OENSC is much higher than that of other competitors. Although the computation time of OENSC-S is slightly longer than that of OENSC, it is also faster than other methods for large-scale datasets. 
Particularly, the traditional batch methods, such as SSC and LRR, struggle with scalability, requiring significantly more time on larger datasets like USPS compared to the proposed online methods. This further emphasizes the advantage of OENSC and OENSC-S in handling large and dynamic datasets effectively.

\section{Conclusion}\label{sec_con}

In this paper, we propose a new subspace clustering framework by integrating both the $\ell_0$ norm and the Frobenius norm, termed as $\ell_0$-ENSC. Theoretically, we prove that the proposed model has the block diagonal property, and establish the relationship between KKT points, P-stable points, and optimal solutions. 
To cope with the challenges posed by online data streaming, we developed an online ADMM algorithm for the $\ell_0$-ENSC model, which incorporates a support-point-based dictionary update strategy. This strategy selectively updates dictionary elements to adapt to new data while maintaining computational costs. We also provide theoretical guarantees for selecting the number of support points, and rigorously prove the convergence of the online algorithm, thereby ensuring its reliability for practical applications. Extensive experiments on several public datasets demonstrate the excellence of our method in terms of clustering performance, robustness, and computational efficiency.

Although the proposed $\ell_0$-ENSC has achieved excellent performance, there are still several areas worth exploring. 
On the one hand, the proposed model involves two hyperparameters, which adds complexity and increases the burden of parameter tuning in practical applications. Therefore, it is necessary to study the adaptive parameter selection mechanism which could further enhance the efficiency and generalization of the algorithm. On the other hand, the computation of support points relies on the CCP algorithm, whose computational efficiency is not satisfactory. Therefore, developing a fast algorithm to find support points with theoretical guarantees, is an important direction for future research.

\section*{Acknowledgment}
This work was supported in part by the 111 Project of China (B16002) and the National Natural Science Foundation of China (12071022).

\section*{Appendix}
\subsection*{Proof of Theorem \ref{th_bdq}}

\begin{proof}
	According to \cite{Lu2019bdr}, we only need to verify that $g(X) = \lambda _1\lVert X \rVert _0+\frac{\lambda _2}{2}\lVert X \rVert _{F}^{2}$ satisfies the following two conditions: (1) $g(P^{\top}XT)=g(X)$ for any permutation matrix $P$; (2) $g(X)\geq g(X^B)$ with $X = \left(\begin{matrix} X_1 & X_2\\ X_3 & X_4 \end{matrix}\right)$ and $X = \left(\begin{matrix} X_1 & 0\\ 0 & X_4 \end{matrix}\right)$, in which the equality holds if and only if $X=X^B$.
	
	For the first part, let $P$ be any permutation matrix, then we have
	\begin{align*}
		g(P^{\top}XP) & =\lambda _1\lVert P^{\top}XP \rVert _0+\frac{\lambda _2}{2}\lVert P^{\top}XP \rVert _{F}^{2}\\
		&= \lambda _1 \sum_{i,j} \mathbb{I} \left[(P^{\top}XP)_{ij}\neq 0\right] +\frac{\lambda _2}{2}\sum_{i,j} \left[ \left( P^{\top}XP\right)_{i,j} \right]^2\\
		&= \lambda _1 \sum_{i,j} \mathbb{I} \left[ P^{\top}(XP)_{ij}\neq 0 \right] +\frac{\lambda _2}{2}\sum_{i,j} \left[ P^{\top}\left( XP\right)_{i,j} \right]^2\\
		&= \lambda _1 \sum_{i,j} \mathbb{I} \left[P^{\top}(X)_{ij}\neq 0 \right] +\frac{\lambda _2}{2}\sum_{i,j} \left[ P^{\top}\left( X\right)_{i,j} \right]^2\\
		&= \lambda _1 \sum_{i,j} \mathbb{I} \left[ X_{ij}\neq 0 \right] +\frac{\lambda _2}{2}\sum_{i,j} \left( X_{i,j}\right) ^2\\
		&= \lambda _1\lVert X \rVert _0+\frac{\lambda _2}{2}\lVert X \rVert _{F}^{2}\\
		&=g(X).
	\end{align*}
	This verifies that $g(X)$ satisfies the first condition. For the second part,  we have
	\begin{align*}
		g(X) 
		= \lambda _1 \sum_{i,j} \mathbb{I} \left(X_{ij}\neq 0\right) +\frac{\lambda _2}{2}\sum_{i,j}X_{i,j}^2 \geq \lambda _1 \sum_{i,j} \mathbb{I} \left(X^B_{ij}\neq 0\right) +\frac{\lambda _2}{2}\sum_{i,j}(X^B_{i,j})^2 = g(X^B).
	\end{align*}
	The equality of the above inequality holds if and only if $X=X^B$, that is $X_2=0$ and $X_3=0$. It then follows from \cite[Theorem 3]{Lu2019bdr} that $X^*$ satisfies the block diagonal property.
\end{proof}

\subsection*{Proof of Proposition \ref{pop_1}}
\begin{proof}
	It follows from \eqref{eq_prox0} that the first equation of \eqref{stationary} can be written as
	\begin{eqnarray*}
		Y^*_{ij} = \left[\operatorname{Prox}_{\alpha\lambda_1\|\cdot\|_0}\left( Y^*+\alpha U^* \right)\right]_{ij} 
		= \begin{cases}
			0, & (Y^*+\alpha U^*)_{ij}\leq \sqrt{2\alpha\lambda_1},\\
			(Y^*+\alpha U^*)_{ij}, & (Y^*+\alpha U^*)_{ij}> \sqrt{2\alpha\lambda_1}.
		\end{cases}
	\end{eqnarray*}
	Further, the above equation is equivalent to
	\begin{eqnarray*}
		0 =   Y^* - \operatorname{Prox}_{\alpha\lambda_1\|\cdot\|_0}\left( Y^*+\alpha U^* \right)
		= \left(\begin{matrix} Y^*_T\\ Y^*_{\bar{T}}\end{matrix}\right) - \left(\begin{matrix} 0 \\ Y^*_{\bar{T}} + \alpha U^*_{\bar{T}} \end{matrix}\right) 
		=  \left(\begin{matrix} Y^*_T \\ \alpha U^*_{\bar{T}} \end{matrix}\right), 
	\end{eqnarray*}
	where $T = \{(i,j): Y^*_{ij}\leq \sqrt{2\alpha\lambda_1}\}$ and $\bar{T}$ is the complement set of $T$. Thus, $Y^*_{ij}=0$ for any $(i,j)\in T$ and $\left| Y^*_{ij} \right| > \sqrt{2\alpha\lambda_1}$  for any $(i,j)\in \bar{T}$, which implies the assertion (a). Then the assertion (b) holds by $Y^*=X^*$. This completes the proof.
\end{proof}

\subsection*{Proof of Theorem \ref{th_nc}}
\begin{proof}
	If $(Y^*,X^*)$ is a local minimizer of \eqref{OENSC1}, then it follows from the equivalence of \eqref{OENSC1} and \eqref{OENSC} that $Y^*$ must be a local minimizer to $\eqref{OENSC}$. By \cite[Theorem 10.1]{rockafellar2009variational}, it has
	\begin{eqnarray*}
		0 \in \nabla h(Y^*) + \lambda_1 \partial\lVert Y^* \rVert _0. 
	\end{eqnarray*}
	Then, by taking $U^* = -\nabla h(X^*)$ with $X^* = Y^*$, it immediately yields that $(Y^*,X^*, U^*)$ is a KKT point. Since $h$ is a strongly convex function, the remaining conclusion follows from \cite[Theorem 3.1]{zhou2023revisiting}. This yields the desired statement.
	
\end{proof}

\subsection*{Proof of Theorem \ref{th_skl}}
\begin{proof}
	By the definition of the proximal operator, it follows from the second equation in \eqref{stationary} that 
	\begin{eqnarray*}
		0 \in \alpha\lambda_1 \partial\lVert Y^* \rVert _0 + Y^* - (Y^* + \alpha U^*).
	\end{eqnarray*}
	This together with $\alpha > 0$ immediately yields that $(Y^*, X^*, U^*)$ is a KKT point. 
	Since $h$ is a strongly convex function, \cite[Theorem 3.2]{zhou2023revisiting} leads to $(Y^*, X^*)$ is also a uniquely local minimizer. The proof is completed.
	
\end{proof}

\subsection*{Proof of Lemma \ref{nonincreasing}}

\begin{proof}
	It follows from \eqref{lag} and \eqref{eq:u} that we have
	\begin{align}\label{eq2}
		\mathcal{L}_{\sigma}(\mathbf{y}_i^{k+1}, \mathbf{x}_i^{k+1}, \mathbf{u}_i^{k+1})- \mathcal{L}_{\sigma}(\mathbf{y}_i^{k+1}, \mathbf{x}_i^{k+1}, \mathbf{u}_i^{k})
		& = -\langle \mathbf{u}_i^{k+1}-\mathbf{u}_i^k, \mathbf{x}_i^{k+1}-\mathbf{y}_i^{k+1}\rangle \nonumber\\
		& = - \frac{1}{\sigma}\|\mathbf{u}_i^{k+1}-\mathbf{u}_i^k\|_2^2. 
	\end{align}
	According to \eqref{eq:x1} and \eqref{eq:u}, we have
	\begin{eqnarray}\label{eq_u0}
		\mathbf{u}_i^{k+1} = -\nabla h(\mathbf{x}_i^{k+1}).
	\end{eqnarray}
	Let $\gamma = \lambda_{min}(D^{\top}D)$, it derives
	\begin{align}\label{eq1}
		\mathcal{L}_{\sigma}(\mathbf{y}_i^{k+1}, \mathbf{x}_i^{k+1}, \mathbf{u}_i^{k+1})- \mathcal{L}_{\sigma}(\mathbf{y}_i^{k+1}, \mathbf{x}_i^{k+1}, \mathbf{u}_i^{k}) 
		& \leq - \frac{1}{\sigma}\| \left(D^{\top}D+\lambda_2 I\right) \left( \mathbf{x}_i^{k+1} - \mathbf{x}_i^k \right)\|_2^2 \nonumber\\
		& \leq - \frac{\gamma+\lambda_2}{\sigma}\|  \mathbf{x}_i^{k+1} - \mathbf{x}_i^k \|_2^2,
	\end{align}

	Recall that $\mathcal{L}_{\sigma}(\mathbf{y}_i^{k}, \mathbf{x}_i, \mathbf{u}_i^{k})$ is strongly convex for $\mathbf{x}_i$ with modulus at least $\gamma + \lambda_2+\sigma$, then it implies
	\begin{eqnarray}\label{eq3}
		\mathcal{L}_{\sigma}(\mathbf{y}_i^{k+1}, \mathbf{x}_i^{k+1}, \mathbf{u}_i^{k}) -\mathcal{L}_{\sigma}(\mathbf{y}_i^{k+1}, \mathbf{x}_i^{k}, \mathbf{u}_i^{k}) \leq -\frac{\gamma + \lambda_2+\sigma}{2}\|\mathbf{x}_i^{k+1}-\mathbf{x}_i^k\|_\textrm{F}^2.
	\end{eqnarray}
	
	Although $\mathcal{L}_{\sigma}(\mathbf{y}_i, \mathbf{x}_i^{k}, \mathbf{u}_i^{k})$ is nonconvex with respect to $\mathbf{y}_i$, $\mathbf{y}_i^{k+1}$ is a (local) minimizer that satisfies
	\begin{eqnarray}\label{eq4}
		\mathcal{L}_{\sigma}(\mathbf{y}_i^{k+1}, \mathbf{x}_i^{k}, \mathbf{u}_i^{k}) - \mathcal{L}_{\sigma}(\mathbf{y}_i^{k}, \mathbf{x}_i^{k}, \mathbf{u}_i^{k})\leq 0.
	\end{eqnarray}
	
	From \eqref{eq1}-\eqref{eq4}, it is not hard to conclude that
	\begin{eqnarray}\label{decrea}
		\mathcal{L}_{\sigma}(\mathbf{y}_i^{k+1}, \mathbf{x}_i^{k+1}, \mathbf{u}_i^{k+1})-\mathcal{L}_{\sigma}(\mathbf{y}_i^{k}, \mathbf{x}_i^{k}, \mathbf{u}_i^{k})\leq - \frac{(2+\sigma)(\gamma+\lambda_2)+\sigma^2}{2\sigma}\|\mathbf{x}_i^{k+1}-\mathbf{x}_i^k\|_2^2.
	\end{eqnarray}
	This yields the desired statement.
\end{proof}

\subsection*{Proof of Lemma \ref{bound}}

\begin{proof}
	
	By the fact that $\mathcal{L}_{\sigma}(\mathbf{y}_i^{k+1}, \mathbf{x}_i^{k+1}, \mathbf{u}_i^{k+1})$ is nonincreasing in Lemma \ref{nonincreasing}, it has
	\begin{eqnarray}\label{eq_bound_1}
		\begin{aligned}
			\mathcal{L}_{\sigma}(\mathbf{y}_i^{k+1}, \mathbf{x}_i^{k+1}, \mathbf{u}_i^{k+1}) &\leq \mathcal{L}_{\sigma}(\mathbf{y}_i^{k+1}, \mathbf{x}_i^{k+1}, \mathbf{u}_i^{k+1}) + \kappa\|\mathbf{x}_i^{k+1}-\mathbf{x}_i^k\|_2^2\\ 
			&\leq \mathcal{L}_{\sigma}(\mathbf{y}_i^{k}, \mathbf{x}_i^{k}, \mathbf{u}_i^{k}) \leq \cdots \\ 
			& \leq \mathcal{L}_{\sigma}(\mathbf{y}_i^{0}, \mathbf{x}_i^{0}, \mathbf{u}_i^{0}),
		\end{aligned}
	\end{eqnarray}
	where the first inequality is derived from \eqref{decrea}. In addition, for $k>1$, the following statements hold
		\begin{align}\label{eq_Lk}
			&\mathcal{L}_{\sigma}(\mathbf{y}_i^{k+1}, \mathbf{x}_i^{k+1}, \mathbf{u}_i^{k+1}) \nonumber\\
			&= h(\mathbf{x}_i^{k+1}) +\lambda _1\lVert \mathbf{y}_i^{k+1} \rVert _0 + \langle \mathbf{u}_i^{k+1}, \mathbf{x}_i^{k+1}-\mathbf{y}_i^{k+1}\rangle + \frac{\sigma}{2}\lVert \mathbf{x}_i^{k+1}-\mathbf{y}_i^{k+1} \rVert _{2}^{2} \nonumber\\
			&= h(\mathbf{x}_i^{k+1}) +\lambda _1\lVert \mathbf{y}_i^{k+1} \rVert _0 + \frac{\sigma}{2}\lVert \mathbf{x}_i^{k+1}-\mathbf{y}_i^{k+1} +  \mathbf{u}_i^{k+1}/\sigma \rVert _{2}^{2} - \frac{1}{2\sigma} \lVert \mathbf{u}_i^{k+1} \rVert _{2}^{2}  \nonumber\\
			&\overset{\eqref{eq_u0}}{=} h(\mathbf{x}_i^{k+1}) +\lambda _1\lVert \mathbf{y}_i^{k+1} \rVert _0 + \frac{\sigma}{2}\lVert \mathbf{x}_i^{k+1}-\mathbf{y}_i^{k+1} -\nabla h(\mathbf{x}_i^{k+1})/\sigma \rVert _{2}^{2} - \frac{1}{2\sigma} \lVert \nabla h(\mathbf{x}_i^{k+1}) \rVert _{2}^{2}  \nonumber\\
			&=  \frac{1}{2} h(\mathbf{x}_i^{k+1}) + \frac{1}{2} \left( h(\mathbf{x}_i^{k+1}) - \frac{1}{2r} \lVert \nabla h(\mathbf{x}_i^{k+1}) \rVert _{2}^{2} \right)  + \left( \frac{1}{4r} - \frac{1}{2\sigma} \right)\lVert \nabla h(\mathbf{x}_i^{k+1}) \rVert _{2}^{2} 
			\nonumber \\
			&~~~ + \lambda _1\lVert \mathbf{y}_i^{k+1} \rVert _0 + \frac{\sigma}{2}\lVert \mathbf{x}_i^{k+1}-\mathbf{y}_i^{k+1} -\nabla h(\mathbf{x}_i^{k+1})/\sigma \rVert _{2}^{2} \nonumber\\			
			& \overset{\eqref{eq_smooth}}{\geq} \frac{1}{2} h(\mathbf{x}_i^{k+1})  + \frac{1}{2} h(\mathbf{x}_i^{k+1} -\frac{1}{r_i}\nabla h_i(\mathbf{x}_i^{k+1})) + \left( \frac{1}{4r} - \frac{1}{2\sigma} \right)\lVert \nabla h(\mathbf{x}_i^{k+1}) \rVert _{2}^{2} 
			\nonumber \\
			&~~~ + \lambda _1\lVert \mathbf{y}_i^{k+1} \rVert _0 + \frac{\sigma}{2}\lVert \mathbf{x}_i^{k+1}-\mathbf{y}_i^{k+1} -\nabla h(\mathbf{x}_i^{k+1})/\sigma \rVert _{2}^{2} \nonumber\\
			& \geq \frac{1}{2} h(\mathbf{x}_i^{k+1})  + \frac{1}{2} \inf_{\mathbf{x}_i} h(\mathbf{x}_i) + \left( \frac{1}{4r} - \frac{1}{2\sigma} \right)\lVert \nabla h(\mathbf{x}_i^{k+1}) \rVert _{2}^{2} 
			\nonumber \\
			&~~~ + \lambda _1\lVert \mathbf{y}_i^{k+1} \rVert _0 + \frac{\sigma}{2}\lVert \mathbf{x}_i^{k+1}-\mathbf{y}_i^{k+1} -\nabla h(\mathbf{x}_i^{k+1})/\sigma \rVert _{2}^{2} \nonumber\\
			& \geq \frac{1}{2} h(\mathbf{x}_i^{k+1})  + \frac{1}{2} \inf_{\mathbf{x}_i} h(\mathbf{x}_i) + \left( \frac{1}{4r} - \frac{1}{2\sigma} \right)\lVert \nabla h(\mathbf{x}_i^{k+1}) \rVert _{2}^{2},
		\end{align}
	where the last inequality follows from the nonnegativity of $\|\cdot\|_0$ and $\|\cdot\|_2$.
	Therefore, substituting \eqref{eq_Lk} into \eqref{eq_bound_1}, we have
	\begin{eqnarray}
			\mathcal{L}_{\sigma}(\mathbf{y}_i^{0}, \mathbf{x}_i^{0}, \mathbf{u}_i^{0}) - \frac{1}{2} \inf_{\mathbf{x}_i} h(\mathbf{x}_i)
			\geq \frac{1}{4}\lVert \mathbf{z}_i-D\mathbf{x}_i^{k+1} \rVert _{2}^{2} +\frac{\lambda_2}{4}\lVert \mathbf{x}_i^{k+1} \rVert _{2}^{2} + \left( \frac{1}{4r} - \frac{1}{2\sigma} \right)\lVert \nabla h(\mathbf{x}_i^{k+1}) \rVert_{2}^{2},  
	\end{eqnarray}
	which immediately yields that $\mathbf{x}_i^{k+1}$ and $ \nabla h(\mathbf{x}_i^{k+1})$  are bounded for  $ \sigma \geq  2 r$. Then, the boundedness of $\mathbf{u}_i^{k+1}$ can be derived from \eqref{eq_u0}. According to \eqref{eq:u}, we have
	\begin{equation}
		\| \mathbf{y}_i^{k+1} \|_2 \leq \|\mathbf{x}_i^{k+1}\|_2 + \frac{1}{\sigma}(\| \mathbf{u}_i^{k+1} \|_2 + \|\mathbf{u}_i^{k}\|_2),
	\end{equation}
	which implies that $\mathbf{y}_i^{k+1}$ is bounded. Overall, the sequence $(\mathbf{y}_i^{k+1}, \mathbf{x}_i^{k+1}, \mathbf{u}_i^{k+1})$ is bounded. 
		
	Now we turn to prove the second part of this lemma. Firstly, it can be derived that the augmented Lagrangian sequence 
	$\mathcal{L}_{\sigma}(\mathbf{y}_i, \mathbf{x}_i, \mathbf{u}_i)$ is bounded below due to the non-negativity of the objective function \eqref{OENSC1} and the boundedness of $\{\mathbf{w}_i^k\}$. This, together with  its nonincreasing given in Lemma \ref{nonincreasing}, demonstrates that $\mathcal{L}_{\sigma}(\mathbf{y}_i, \mathbf{x}_i, \mathbf{u}_i)$ must converge. Thus, $\lim_{k\rightarrow \infty} \|\mathbf{x}_i^{k+1}-\mathbf{x}_i^k\|_2^2=0$ can be derived by taking limit on both sides of \eqref{decrea}. Similarly, from \eqref{eq2}, it obtains $\lim_{k\rightarrow \infty} \|\mathbf{u}_i^{k+1}-\mathbf{u}_i^k\|_\textrm{F} = 0$.  It follows from \eqref{eq:u} that
	\begin{eqnarray*}
		\mathbf{y}_i^{k+1}-\mathbf{y}_i^k = \mathbf{x}_i^{k+1}-\mathbf{x}_i^k  - \frac{1}{\sigma}(\mathbf{u}_i^{k+1} - \mathbf{u}_i^{k}) + \frac{1}{\sigma}(\mathbf{u}_i^{k}-\mathbf{u}_i^{k-1}),
	\end{eqnarray*}
	which implies that
	\begin{eqnarray*}
		\|\mathbf{y}_i^{k+1}-\mathbf{y}_i^k\|_2 \leq \|\mathbf{x}_i^{k+1}-\mathbf{x}_i^k\|_2  + \frac{1}{\sigma}\|\mathbf{u}_i^{k+1} - \mathbf{u}_i^{k}\|_2 
		+ \frac{1}{\sigma}\|\mathbf{u}_i^{k}-\mathbf{u}_i^{k-1}\|_2.
	\end{eqnarray*}
	By taking limits on both sides of the above equation, we have
	$\lim_{k\rightarrow \infty}\|\mathbf{y}_i^{k+1}-\mathbf{y}_i^k\|_2 = 0$.  This completes this proof.	
\end{proof}

\subsection*{Proof of Theorem \ref{th_con}}
\begin{proof}
	The proof of this theorem proceeds in three steps. Firstly, we prove that $\{\mathbf{w}_i^k\}$ has a limit point. It follows from Lemma \ref{bound} that the generated sequence $\{\mathbf{w}_i^k\}$ is bounded and the augmented Lagrangian sequence $\mathcal{L}_{\sigma}(\mathbf{y}_i^k, \mathbf{x}_i^k, \mathbf{u}_i^k)$ is convergent. Therefore, $\{\mathbf{w}_i^k\}$ has at least one limit point. Furthermore, for any $\{\mathbf{w}_i^k\}$ in the set of all limit points, there exists a convergent subsequence $\{\mathbf{w}_i^{k_j}\}$ which satisfies $\lim\limits_{k_j\rightarrow \infty} \mathbf{w}_i^{k_j} = \mathbf{w}_i^*.$ 
	
	Secondly, we show that any limit points of $\{\mathbf{w}_i\}$ is a P-stationary point of problem \eqref{OENSC1}. 
	By taking the limit on both sides of the update formula for the Lagrangian multiplier, it can be reached that $\mathbf{x}_i^* - \mathbf{y}_i^*=0$, which demonstrates the feasibility of $\mathbf{w}_i^*$. 
	Then, the first equation in \eqref{stationary} can be obtained directly from the $\mathbf{x}_i$-subproblem. 
	Furthermore, $\mathbf{y}_i^{k_j}  = \operatorname{Prox}_{\frac{\lambda_1}{\sigma}\|\cdot\|_0}(\mathbf{x}_i^{k_j} + \mathbf{u}_i^{k_j}/\sigma)$ can be derived by taking limits on both	sides of \eqref{eq:y}. Thus, $\{\mathbf{w}_i^*\}$ is a stationary point of problem \eqref{OENSC1}. 
	
	Finally, $\{\mathbf{w}_i^*\}$ is a unquely local optimal solution of \eqref{OENSC1} which can be derived from Theorem \ref{th_skl} immediately. 
	In summary, this yields the desired conclusion.
\end{proof}

\subsection*{Proof of Lemma \ref{le_convex}}
\begin{proof}
	It's obvious that $Q(\mathbf{d}^{'}|\mathbf{d}^{'}) =\|\mathbf{d}^{'}\|_2$. For $\mathbf{d}\neq \mathbf{d}^{'}$, it has
		\begin{eqnarray*}
			(\| \mathbf{d} \|_2 - \| \mathbf{d}^{'} \|_2 )^2 = \| \mathbf{d} \|_2^2 + \| \mathbf{d}^{'} \|_2^2 - 2 \| \mathbf{d} \|_2 \| \mathbf{d}^{'} \|_2 \geq 0
		\end{eqnarray*}
	which implies $\| \mathbf{d} \|_2^2 + \| \mathbf{d}^{'} \|_2^2 \geq 2 \| \mathbf{d} \|_2 \| \mathbf{d}^{'} \|_2$. Thus,
	\begin{eqnarray*}
		Q(\mathbf{d}|\mathbf{d}^{'})  = \frac{\| \mathbf{d} \|_2^2}{2 \| \mathbf{d}^{'} \|_2} + \frac{\| \mathbf{d}^{'} \|_2}{2} = \frac{\| \mathbf{d} \|_2^2 + \| \mathbf{d}^{'} \|_2^2}{2 \| \mathbf{d}^{'} \|_2} \geq \frac{2 \| \mathbf{d}^{'} \|_2 \| \mathbf{d} \|_2}{2 \| \mathbf{d}^{'} \|_2} = \| \mathbf{d} \|_2.
	\end{eqnarray*}
	This yields the desired statement.
\end{proof}

\bibliography{mybibfile}

\end{document}